\documentclass[12pt]{article}
\usepackage{latexsym,amssymb,amsmath,a4wide}

\newtheorem{proposition}{Proposition}
\newtheorem{theorem}{Theorem}

\newcommand{\Aut}{{\rm Aut}}

\newcommand{\veps}{\varepsilon}

\newcommand{\ord}{{\rm ord}}

\begin{document}

\title{\vspace{-3cm} Orientably-regular embeddings of complete multigraphs}
\author{}
\date{}
\maketitle

\begin{center}
\vspace{-1.3cm}

{\large \v{S}tefan Gy\"urki} \\
\vspace{1mm} {\small Slovak University of Technology, Bratislava, Slovakia}\\

\vspace{4mm}

{\large So\v{n}a Pavl\'ikov\'a} \\
\vspace{1mm} {\small A. Dub\v{c}ek University, Tren\v{c}\'in, Slovakia}\\

\vspace{4mm}

{\large Jozef \v Sir\'a\v n} \\
\vspace{1mm} {\small Slovak University of Technology, Bratislava, Slovakia}

\vspace{4mm}

\end{center}

\begin{abstract}
An embedding of a graph on an orientable surface is {\em orientably-regular} (or {\em rotary}, in an equivalent terminology) if the group of orientation-preserving automorphisms of the embedding is transitive (and hence regular) on incident vertex-edge pairs of the graph. A classification of orientably-regular embeddings of complete graphs was obtained by L. D. James and G. A. Jones [in "Regular orientable imbeddings of complete graphs", {\em J. Combinatorial Theory Ser. B 39 (1985), 353--367}], pointing out interesting connections to finite fields and Frobenius groups. By a combination of graph-theoretic methods and tools from combinatorial group theory we extend results of James and Jones to classification of orientably-regular embeddings of complete multigraphs with arbitrary edge-multiplicity.
\end{abstract}

\section{Introduction}\label{sec:intro}

Cellular embeddings of connected graphs on compact orientable surfaces, or, simply, {\em maps}, provide important links  between discrete mathematics, low-dimensional topology and Riemann surfaces, extending to group theory and Galois theory if the maps are `highly symmetric'. To be more specific about the latter, an {\em automorphism of a map} is an automorphism of its {\em underlying graph} that preserves the cell structure of a map and orientation of its {\em carrier surface}; a map is {\em orientably-regular} if the group of all its automorphisms acts regularly on incident vertex-edge pairs of its underlying graph. Orientable regularity is known to be a formal way of expressing the  `highest level of symmetry' a map can have with respect to orientation-preserving automorphisms, cf. \cite{JoSi}.
\smallskip

Classification of orientably-regular maps, both per se as well as because applications to symmetric Riemann surfaces, turns out to be an important challenge and has been pursued for more than a century. Since the concept of an orientably-regular map encompasses an interplay of three kinds of objects (graphs, surfaces, and  groups), classification of such maps broadly follows three directions, depending on the object of primary interest. Leaving aside a number of deep results obtained in classification of orientably-regular maps with a given automorphism group, or on a given surface, we confine ourselves only to a brief survey on classification results for orientably-regular maps with a given underlying graph, focusing on important infinite families of graphs.
\smallskip

The four highlights in this area of research are classifications of orientably-regular embeddings of complete graphs \cite{JaJo}, complete bipartite graphs \cite{Jones}, complete multipartite graphs with at least three parts \cite{DuZh}, and cubes of arbitrary finite dimension \cite{Cat+}. In all four cases the results are highly non-trivial, despite targeting well-studied families of graphs. Indeed, one may naively ask about possible difficulties in sorting out orientably-regular embeddings of complete graphs in which everything should be understood. The paper \cite{JaJo}, which is historically first in classification of orientably-regular maps with a given underlying graph, presents a surprising answer which includes deep connections of the problem with finite fields and Frobenius groups. By its main result, a complete graph on $r$ vertices admits an orientably-regular embedding if and only if $r=p^k$ for some prime $p$, in which case there are exactly $\varphi(r-1)/k$ pairwise non-isomorphic orientably-regular embeddings (where $\varphi$ is Euler's function).
\smallskip

For completeness, although a general characterisation of graphs underlying orientably-regular maps  in terms of existence of a suitable subgroup of the automorphism group of the graph was given in \cite{GNSS}, the four major results mentioned above illustrate a long way from this abstract characterisation to results for particular families of graphs.
\smallskip

The aim of this paper is to extend the classification of orientably-regular embeddings of complete graphs $K_r$ to complete {\em multigraphs} $K_r^{(t)}$ of order $r$ in which every pair of vertices is joined by $t$ parallel edges for some $t\ge 1$; we also speak of a complete graph with {\em edge-multiplicity} $t$. A classification of such embeddings in the range $2\le r\le 5$ for every $t\ge 1$ was obtained in \cite{Gar} as part of a project to determine all orientably-regular maps of order at most $5$; the cases $r=2$ and $r=3$ were later revisited in \cite{NeSk} and \cite{H+} in different contexts. But beyond the above range of $r$, classification remained completely open. In this place we offer the following very brief excerpt from our main findings.

\begin{theorem}\label{thm:intro}
A complete multigraph $K_r^{(t)}$ with $r\ge 5$ and edge-multiplicity $t\ge 1$ admits an orientably-regular embedding if and only if $r$ is a prime power, with $t$ odd if $r$ is odd. Moreover, if $r=p^k$ for some prime $p$ and $k\ge 1$, the number of mutually non-isomorphic orientably-regular embeddings of $K_r^{(t)}$ is $\varphi(r-1)/k$ for every odd $t$ and arbitrary $r\ge 5$, and $2\varphi(r-1)/k$ for every even $t\ge 2$ and every even $r\ge 8$.
\end{theorem}

Our proofs rely on a combination of graph-theoretic arguments and tools from combinatorial group theory. We actually give two refined versions of the above theorem: one in terms of central cyclic extensions of automorphism groups of regular embeddings of simple complete graphs (that is, those with edge-multiplicity $1$), and a more detailed one in terms of group presentations. For the latter we construct a suitable associated Cayley graph and apply the Reidemeister-Schreier procedure, which is based on converting closed walks in the associated graph to relators in group presentations. The methods apply to every $r\ge 3$, furnishing also an explanation (via Schur multiplier) of why the value of $r=4$ is a remarkable exception, as it also follows from \cite{Gar}.
\smallskip

The outlined approach determines the structure of the paper. After a presentation of basic facts on orientably-regular maps (section \ref{sec:basic}) we review specifications of the theory to underlying multigraphs with arbitrary edge-multiplicity (section \ref{sec:multi}). This has been kind of folklore but we are not aware of a synthetic approach we offer here, linking such embeddings with a combinatorial description of {\em presentations} of the corresponding automorphism groups. In section \ref{sec:simple} we refine the known classification of orientably-regular embeddings of simple complete graphs of prime power order by deriving presentations of their automorphism groups, which we then extend to corresponding complete multigraphs of arbitrary edge-multiplicity for prime powers larger than $4$ in section \ref{sec:multi-compl}. The `singular values' $2$, $3$ and $4$ are then handled separately in section \ref{sec:sing}, and we wrap up in section \ref{sec:conc} by a discussion of chirality, self-complementarity, Wilson operators, and Cayley structure of the maps.
\smallskip

One may also wonder why one should be interested in investigation of orientably-regular embeddings of {\em multigraphs} at all. The perhaps best answer is provided by an inspection of a computationally generated list of presentations of automorphism groups of all orientably-regular maps on surfaces up to genus $300$ \cite{Co-list}, in which only $11$ per cent are such that both the map and its dual have a simple underlying graph. It appears that, at least for genera up to $300$, orientably-regular maps with multiple edges are vastly predominant.

\section{Orientably-regular maps: a brief review of theory}\label{sec:basic}

Let us begin by recalling that in the theory of maps, automorphisms of graphs and maps are usually viewed as permutations of incident vertex-edge pairs, commonly known as {\em arcs}. An automorphism of a map is thus a permutation of the set of arcs of its underlying graph which is both an automorphism of the graph and also preserves face boundary walks of the map, including their orientations inherited from an orientation of the carrier surface of the map. The group of all such automorphisms of a map ${\cal M}$ will be denoted $\Aut^+{\cal M}$, the `plus' superscript indicating that the group comprises only orientation-preserving automorphisms; we will address their orientation-reversing counterparts in section \ref{sec:conc}.
\smallskip

It is an elementary fact in the theory of maps that, for a map ${\cal M}$ as above, the group $\Aut^+{\cal M}$ acts freely on the set of arcs of the underlying graph of ${\cal M}$. In the case the action is also transitive, and hence regular, we say that ${\cal M}$ is an {\em orientably-regular} map. Equivalently, in an orientably-regular map ${\cal M}$, for each ordered pair of arcs there exists an automorphism in $\Aut^+{\cal M}$ sending the first arc of the pair onto the second.
\smallskip

Let now ${\cal M}$ be an orientably-regular map, formed by a cellular embedding of its underlying graph $\Gamma$ on a compact orientable surface ${\cal S}$, oriented anti-clockwise, and let $a$ be an arbitrary but fixed arc of $\Gamma$. By orientable regularity, there is an automorphism in $\Aut^+{\cal M}$, which we denote $y$, mapping the arc $a$ to the counterclockwise next arc $a'$ of ${\cal M}$ incident to the same vertex as $a$. Similarly, orientable regularity implies existence of another automorphism in $\Aut^+{\cal M}$, denoted $x$ in what follows, which sends the arc $a$ onto the anti-clockwise next arc on the face of ${\cal M}$ containing both $a$ and $a'$. It may be checked that $a(xy)^2=a$, and invoking orientable regularity again this means that $(xy)^2=1$ in $\Aut^+{\cal M}$. Because of the facts that the automorphisms $x$, $y$ and $xy$ act on ${\cal M}$ as local rotations around centre of a face, a vertex, and centre of an edge (respectively), in an alternative terminology such maps are also known as {\em rotary} \cite{Wi}.
\smallskip

Continuing with our orientably-regular map ${\cal M}$, connectedness of $\Gamma$ (and hence of the carrier surface of the map) implies that the entire group $\Aut^+{\cal M}$ is {\em generated} by the two permutations $x$ and $y$. Further, and again by orientable regularity of the map, all its faces must be bounded by boundary walks of the same length, say, $m$, and all its vertices must have the same valency, say, $n$; using the traditional Schl\"afli symbol the map is then said to be of {\em type} $\{m,n\}$. The connection of the type to the generators $x$ and $y$ of $\Aut^+{\cal M}$ is simply that entries of the type are orders of the generators, that is, $m=\ord(x)$ and $n=\ord(y)$. The group $\Aut^+{\cal M}$ thus admits a presentation of the form $\langle x,y\ |\ x^m,\,y^n,\,(xy)^2,\,\ldots\rangle$. This fact has a number of consequences, out of which we outline only two in what follows.
\smallskip

First, as any structure preserved by a group acting regularly on building blocks of the structure, the map ${\cal M}$ is completely determined by the group $G=\Aut^+{\cal M}=\langle x,y\rangle$ by identifying arcs with elements of $G$ and faces, vertices and edges of the map with left cosets of the subgroups $\langle x\rangle$, $\langle y\rangle$ and $\langle xy\rangle$, respectively, with incidence defined by containments of elements in cosets  and non-empty intersection of cosets, and with automorphisms acting by right multiplication. This way, orientably-regular maps of type $\{m,n\}$ may be identified with finite groups generated by two elements, of order $m$ and $n$, with product of order $2$.
\smallskip

Second, bearing in mind such an identification, a group $G$ given by a presentation $\langle x,y\ |\ x^m,\,y^n,\, (xy)^2,\,\ldots\rangle$ is a smooth quotient of the well-known {\em ordinary $(m,n,2)$-triangle group} $\Delta^+(m,n,2) = \langle X,Y\ |\ X^m,\,Y^n,\,(XY)^2\rangle$. The latter is the group of orientation-preserving automorphisms of a {\em universal tessellation} $U(m,n)$ of a simply-connected surface by $m$-gons, $n$ meeting at each vertex; the surface is a sphere, a Euclidean plane or a hyperbolic plane, depending on whether $1/m+1/n-1/2$ is positive, zero, or negative. The smooth epimorphism $\Delta^+(m,n,2)\to G$ given by $X\mapsto x$ and $Y\mapsto y$ translates to a smooth covering of the orientably-regular map ${\cal M}$ determined by $G$ by the tessellation $U(m,n)$, which can be used to `pull down' the complex structure from a simply-connected surface carrying $U(m,n)$ onto the carrier surface of ${\cal M}$, making it a Riemann surface.
\smallskip

For further facts related to the theory of orientably-regular maps and their relations to two-generated groups, Riemann surfaces and Galois theory we refer to the seminal paper \cite{JoSi} laying foundations of the theory of maps on orientable surfaces, and to a more recent survey \cite{Si-surv}. The latter also addresses {\em regular maps}, which, on orientable surfaces, arise if an orientably-regular map admits also orientation-reversing automorphisms (and the concept extends also to non-orientable surfaces); these will be briefly mentioned in section \ref{sec:conc}.
\smallskip

\section{Orientably-regular maps with multiple edges}\label{sec:multi}

As alluded to in the introduction, in this section we collect known facts on relations between orientably-regular maps with a simple underlying graph and orientably-regular embeddings of the corresponding multigraph, obtained by `inflating' adjacencies in the original graph to any multiplicity $t\ge 1$. We reiterate that most of the material has been known, but we are not aware of a synthesis combining ideas from graph covers and combinatorial group theory to derive a presentation of the automorphism group of the `inflated' map terms of the automorphism group of the original map.
\smallskip

Let $\cal M$ be an orientably-regular map of type $\{m,n\}$, with its group $G=\Aut^+{\cal M} = \langle x,y\ |\ x^m,y^n,(xy)^2,\ldots \rangle$ of all orientation-preserving automorphisms. Assume that the underlying graph $\Gamma$ of this map contains a pair of distinct vertices joined by exactly $t$ parallel edges for some $t\ge 2$. Since $G$ acts on the graph arc-transitively, {\em every} pair of adjacent vertices in $\Gamma$ is joined by $t$ parallel edges. Thus, $\Gamma$ has {\em edge-multiplicity} $t$ and valency $n=dt$ for some integer $d$, so that every vertex in $\Gamma$ has $d=n/t$ neighbours. Assuming that $y$ is an $n$-fold rotation of $\cal M$ about a chosen vertex of $\Gamma$, the subgroup $\langle y\rangle < G$ of order $n$ acts cyclically on the $d$ neighbours of that vertex.
\smallskip

It follows that $y^d$ fixes not only the chosen vertex but also all its neighbours, and as ${\cal M}$ is orientably-regular, $y^d$ fixes every vertex of the map and so $N=\langle y^d\rangle$ is a normal subgroup of $G$. The orientably-regular quotient map ${\cal M}/N$ then has a simple underlying graph, namely, the quotient $\Gamma/N$, of valency $d=n/t$. The natural projection ${\cal M}\to {\cal M}/N$ induces a branched covering of the two maps, with branching index $t$ at every vertex.
\smallskip

One can turn around this line of thought and begin with an orientably-regular map determined by a group $G=\langle x,y\rangle$ as above and {\em assume} that $N=\langle y^d\rangle$ is a normal subgroup of $G$ for some divisor $d$ of $n$. The subgroup $N$ again fixes every vertex of $\cal M$, and the quotient map ${\cal M}/N$ will be orientably-regular, with the same number of vertices as $\cal M$, and of valency $d$. The underlying graph of this map is simple if and only if $N$ is the largest subgroup of $\langle y\rangle$ normal in $G$, that is, if and only if $N$ is the {\em core} of $\langle y\rangle$ in $G$.
\smallskip

Returning to our original setting with $n=dt$ and edge-multiplicity $t$, note that conjugation of $y^d$ by both $x$ and $xy$ results in the same power of $y^d$, say,  $y^{df}$, which gives rise to a relator of the form $R=x^{-1}y^d xy^{-df}$ in the presentation of $G$. But because $xy$ is an involution, it follows that $f^2\equiv 1$ mod $t$, with a consequence that $[x^2,y^d]=1$. Moreover, if $\Gamma$ (and, equivalently, $\Gamma/N$) is not bipartite, then $f\equiv 1$ mod $t$. Indeed, $\Gamma$ then contains a cycle of odd length $j$, which corresponds to a relator $W$ in $G$ of the form $W=y^{\ell_1}xy^{\ell_2}x\ldots y^{\ell_j}x$. Applying $R$ in the equivalent form $y^dx = xy^{df}$ successively $j$ times to the product $y^dW$ to `move the term $y^d$ across $W$ from the left to the right' results in  $y^dW= Wy^{df^j}$, and since $W=1$ and $j$ is odd, it follows that $f\equiv 1$ mod $t$. In particular, if $\Gamma$ is not bipartite, then the relator $R$ reduces to the commutator $[x,y^d]=x^{-1}y^{-d}xy^d$.
\smallskip

Letting $H=G/N$ and $u=xN$, $v=yN$, the group $H\cong \Aut^+{\cal M}'$ of the quotient map ${\cal M}'={\cal M}/N$ is generated by $u$ and $v$ and admits a presentation of the form
\begin{equation}\label{eq:H} H=\langle u,v\ |\ u^\ell,v^d,(uv)^2,{\rm Rels}\rangle \end{equation}
for some divisor $\ell$ of $m$ and a finite collection $\rm Rels$ of relators as words in $u$ and $v$.
\smallskip

We will now be interested in the reverse direction. Namely, assume that ${\cal M}'$ is an orientably-regular map of type $\{\ell,d\}$ with a simple underlying graph, given by means of a presentation of its group $H=\Aut^+{\cal M}'$ as in \eqref{eq:H}. For any given integer $t\ge 2$, we would like to describe all orientably-regular maps ${\cal M}$ of type $\{m,n\}$ for $n=dt$ (where $m$ is some multiple of $\ell$) such that the underlying graph of ${\cal M}$ is obtained by replacing every edge of the underlying graph of ${\cal M}'$ by $t$ parallel edges. Each such orientably-regular map ${\cal M}$ will be briefly called a $t$-{\em inflation} of ${\cal M}'$, and in this terminology the task is to describe all $t$-inflations of a given orientably-regular map with a simple underlying graph.
\smallskip

The outlined theory implies that such a task reduces to describing suitable normal extensions of $H$ by a cyclic group of order $t$; the procedure is well known and we briefly describe its essence. Let $\langle u,v \ |\ u^\ell, v^d, (uv)^2,{\rm Rels}\rangle$ be a presentation of $H=\Aut^+{\cal M}'$ for an orientably-regular map ${\cal M}'$ of type $\{\ell,d\}$ as in \eqref{eq:H}. Let ${\cal M}$ be a $t$-inflation of ${\cal M}'$, of type $\{m,n\}$ for $n=dt$ and $m$ a multiple of $\ell$, with $G=\Aut^+{\cal M} = \langle x,y\ |\ x^m,y^n,(xy)^2,\ldots\rangle$. We saw that the two maps and the two groups are related by ${\cal M}'= {\cal M}/N$ and $H=G/N$ for the normal subgroup $N=\langle y^d\rangle$ of order $t$. Now, for every word $w(u,v)$ in $\{u^\ell\}\cup\rm Rels$, let $w(x,y)$ be the corresponding word in terms of $x$ and $y$ obtained from $w(u,v)$ by substituting $x$ and $y$ for every occurrence of $u$ and $v$ in $w(u,v)$, respectively. Since $w(x,y)N = w(xN,yN)= w(u,v)=1$, it follows that $w(x,y)\in N$, and so $w(x,y)=(y^d)^{e_w}$ for some integer $e_w$ depending on $w=w(u,v)\in \{u^\ell\}\cup{\rm Rels}$.
\smallskip

By \cite[Ch. 10, Proposition 1]{John} adapted to maps setting, if an orientably-regular map ${\cal M}$ is a $t$-inflation of an orientably-regular ${\cal M}'$ of type $\{\ell,d\}$ with $H=\Aut{\cal M}'$ presented as in \eqref{eq:H}, then the group $G=\Aut^+{\cal M}$ admits a presentation of the form
\begin{equation}\label{eq:lift}
G=\langle x,y\ |\ y^n,(xy)^2,x^{-1}y^d xy^{-df},\{w(x,y)y^{-de_w};\ w=w(u,v)\in \{u^\ell\}\cup{\rm Rels}\}\rangle
\end{equation}
for some integers $e_w$ for $w\in \{u^\ell\}\cup{\rm Rels}$ and $f$ such that $f^2\equiv 1$ mod $t$, with $n=dt$.
\smallskip

Shifting focus to inflations of {\em complete maps}, let ${\cal M}$ be an orientably-regular embedding of a complete graph $K_r^{(t)}$ of order $r\ge 2$ with edge-multiplicity $t\ge 1$, and let $G= \Aut^+{\cal M} = \langle x,y\ |\ x^m,y^n, (xy)^2, \ldots\rangle$. In this situation one has $n=dt$ for $d=r-1$, with  $N=\langle y^{d} \rangle \triangleleft G$, and the quotient map determined by the group $H=G/N = \langle u,v\ |\ u^\ell,v^d,(uv)^2, ...,\rangle$ for $u=xN$ and $v=yN$ is an orientably-regular embedding of the (simple) complete graph $K_r$. Orientably-regular embeddings of the latter have been classified in \cite{JaJo}; it turns out that their order is necessarily some prime power $q$ and they are of type $\{\ell,d\}$, where $d=q-1$ and the values of $\ell$ are as follows: $\ell=q$ if $q\in\{2,3\}$, $\ell=q-1$ if $q$ is even but $q\ne 2$ or $q\equiv 1$ mod $4$, and $\ell=(q-1)/2$ if $q\equiv 3$ mod $4$ but $q\ne 3$. Moreover, the group of orientation-preserving automorphisms of such a map turns out to be isomorphic to ${\rm AGL}(1,q)$, the affine general one-dimensional groups over finite fields of order $q$.
\smallskip

Our goal is to classify orientably-regular embeddings ${\cal M} = {\cal M}(q,t)$ of complete graphs of order $q$ with edge-multiplicity $t$ for an arbitrary given prime power $q$ and a positive integer $t$. If such a map is of type $\{m,n\}$, with $n=dt=(q-1)t$, it would be natural to try to identify the maps in terms of their groups $G=\Aut^+{\cal M}(q,t) = \langle x,y\ |\ x^m,y^n, (xy)^2,...\rangle$ by lifting presentations of their quotients $G/N$ for $N=\langle y^d\rangle$ as given by \eqref{eq:lift}. Since $K_q$ is not bipartite for $q>2$, in \eqref{eq:lift} one then has $f=1$, which means that $N$ is a central subgroup of $G$, making the extension of the quotient $G/N$ {\em central} for $q>2$.
\smallskip

The process of lifting, however, assumes knowledge of {\em presentations} of the quotient automorphism groups $H_q= \langle u,v\ |\ u^\ell,v^d, (uv)^2, ...,\rangle\cong {\rm AGL}(1,q)$ of orientably-regular embeddings of the corresponding {\em simple} complete graphs, which appear to be known only for $q\le 5$ from \cite{Gar} as part of a classification of regular maps with underlying graphs of order at most $5$. It should be emphasised that the classification of orientably-regular embeddings of simple complete graphs of \cite{JaJo}, although deep and taken advantage of a number of times here, is {\em not} in terms of group presentations, and hence cannot be directly used for lifting.
\smallskip

Because of this state-of-the-art, in the next section we will use the classification of \cite{JaJo} to determine all the orientably-regular maps with simple complete underlying graphs {\em in terms of presentations of their automorphism groups}.

\section{Presentations of automorphism groups of \\ orientably-regular simple complete maps}\label{sec:simple}

Let $q=p^k$ be an arbitrary prime power and let $F_q$ be a finite field of order $q$. The group $H_q=\langle u,v\rangle \cong {\rm AGL}(1,q)$ is a split extension of its additive (and characteristic) subgroup $F_q^+ \cong (F_p^+)^k$ of $F_q$ by its multiplicative group $F^\times_q$; the subgroup $F_q^+$ is, at the same time, the commutator subgroup of $H$. Let $\mu=\mu(\lambda)=\lambda^k+a_{k-1}\lambda^{k-1} +\ldots +a_1\lambda +a_0$ be a monic primitive polynomial for the field $F_q$ with coefficients in its prime field $F_p$, that is, a monic irreducible polynomial with roots $\{\xi^{p^i};\ 0\le i\le k-1\}$ for some primitive element $\xi\in F_q^\times$ (or, still equivalently, an irreducible factor of the corresponding cyclotomic polynomial over $F_q$). There are $\varphi(q-1)/k$ such primitive polynomials, and it follows from \cite{JaJo} and also from arguments in \cite{JLSW} that these polynomials are in a one-to-one correspondence with isomorphism classes of orientably-regular embeddings of the graph $K_q$.
\smallskip

Recall that there is a standard way to associate a $k\times k$ matrix over $F_p$ (equivalently, a linear transformation of $(F_p^+)^k$ regarded as a $k$-dimensional vector space over $F_p$) the minimal polynomial of which is equal to the polynomial $\mu$ chosen above. For our purposes it will be of advantage to consider the transpose of the standard companion matrix and associate with $\mu$ the matrix $A$ defined by
\begin{equation}\label{eq:m}
A=\left(\begin{array}{ccccccc}
0 & 0 & 0 & ... & 0 & 0 & -a_0 \\
1 & 0 & 0 & ... & 0 & 0 & -a_1 \\
0 & 1 & 0 & ... & 0 & 0 & -a_2 \\
& ... &   & ... &   &  & ...  \\
0 & 0 & 0 & ... & 1 & 0 & -a_{k-2} \\
0 & 0 & 0 &... & ... & 1 & -a_{k-1}
\end{array}\right)
\end{equation}
with ${\rm det}(\lambda I-A)=\mu(\lambda)$. Primitivity of $\mu$ implies that $A$ has order $q-1$ and every element of $H_q\cong F_q^+\rtimes F_q^{\times}\cong {\rm AGL}(1,q)$ has the form $(g,A^i)$ for some $i\in \{0,1,\ldots,q-2\}$ and $g\in (F_p^+)^k$, regarded as a $k$-dimensional column vector, with multiplication in the semidirect product given by $(g,A^i) \cdot(h,A^j) = (g+A^ih,A^{i+j})$. (Note that this description of a semidirect product harmonises with multiplication of column vectors by matrices from the left and with using the transpose of a companion matrix.) In what follows we will denote by $\mathbf{0}$ the all-zero column vector of dimension $k$, and for $1\le i\le k$ we will use the standard notation $\mathbf{e}_i\in (F_p^+)^k$ for the $i$-th unit column vector of dimension $k$. Also, let $\mathbf{a}\in (F_p^+)^k$ denote the column vector forming the last column of $A$.
\smallskip

With this notation we are in position to state a further consequence of \cite{JaJo}, which is that for the generators $u$ and $v$ of the group $H_q\cong {\rm Aut}^+{\cal M}$ one may take
\begin{equation}\label{eq:uv} u = (\mathbf{e}_1,-A^{-1})\ \ \mbox{and} \ \ v=(\mathbf{0},A), \ \ \mbox{with} \ \ uv= (\mathbf{e}_1,-I)  \end{equation}
and routine calculations (see \cite{JLSW} for more details) imply that the order of $u$ is $q-1$ if $q$ is a power of $2$ or $q\equiv 1$ mod $4$, and $(q-1)/2$ if $q\equiv 3$ mod $4$; the orders of $v$ and $uv$ are $q-1$ and $2$, respectively. In addition, one may verify that
\begin{equation}\label{eq:Ae}
A\mathbf{e}_i=\mathbf{e}_{i+1}\ \ \ (1\le i\le k-1) \ \ \mbox{and} \ \ A\mathbf{e}_k=\mathbf{a}
\end{equation}
and, consequently,
\begin{equation}\label{eq:AeI}
A^{i-1}(A-I)\mathbf{e}_1 = \mathbf{e}_{i+1}-\mathbf{e}_i \ \ (1\le i\le k-1) \ \ \mbox{and} \ \ A^{k-1}(A-I)\mathbf{e}_1 = \mathbf{a} - \mathbf{e}_k \ .
\end{equation}
The way $u$ and $v$ have been introduced together with equations \eqref{eq:AeI} imply that $u^{-1}v^{-1}uv=[u,v] = ((A-I)\mathbf{e}_1,I) = (\mathbf{e}_2-\mathbf{e}_1,I)$, and a similar calculation shows that
\begin{equation}\label{eq:[u,v]}
v^{i-1}[u,v]v^{1-i} = (A^{i-1}(A-I)\mathbf{e}_1,I) = (\mathbf{e}_{i+1}-\mathbf{e}_i,I), \ \ 1\le i\le k-1, \ \
\end{equation}
together with
\begin{equation}\label{eq:[u,v]k}
v^{k-1}[u,v]v^{1-k} =  (A^{k-1}(A-I)\mathbf{e}_1,I) = (\mathbf{a}-\mathbf{e}_k,I)\ .
\end{equation}
Let $v_q=1$ if $q$ is even, and $v_q=v^{(q-1)/2}$ if $q$ is odd. We now prove by induction that
\begin{equation}\label{eq:ind}
(\mathbf{e}_i,I) = v^{i-1}(uv)v^{1-i}v_q\ \ \mbox{for every} \ \ i\in\{1,2,\ldots,k\}\ .
\end{equation}
For even $q$ one has \eqref{eq:ind} for $i=1$ trivially, and for odd $q$, from $v^{(q-1)/2}=(\mathbf{0},-I)$ it follows that $(uv)v_q = (\mathbf{e}_1,I)$, which is again \eqref{eq:ind} for $i=1$. The induction hypothesis, that is, $v^{i-1}(uv)v^{1-i}v_q = (\mathbf{e}_i,I)$ for $i\ge 1$ and $i\le k-1$, together with \eqref{eq:[u,v]} implies that
\begin{equation*} (\mathbf{e}_{i+1},I) = (\mathbf{e}_{i+1}-\mathbf{e}_i,I)\cdot (\mathbf{e}_i,I)
= v^{i-1}[u,v]v^{1-i}\cdot v^{i-1}(uv)v^{1-i}v_q
\end{equation*}
and the last expression is easily seen to simplify to $v^i(uv)v^{-i}v_q$, completing the induction step. By almost the same token, combining \eqref{eq:[u,v]k} with the case $i=k$ of \eqref{eq:ind} gives
\begin{equation}\label{eq:a-k}
(\mathbf{a},I) = (\mathbf{a}-\mathbf{e}_k,I)\cdot (\mathbf{e}_k,I) = v^{k-1}[u,v]v^{1-k}\cdot v^{k-1}(uv)v^{1-k}v_q = v^k(uv)v^{-k}v_q\ .
\end{equation}
Observe that \eqref{eq:ind} also implies that the elements $v^{i-1}(uv)v^{1-i}v_q$ for $i\in\{1,2,\ldots,k\}$ commute with each other; this feature can be expressed in a unified form by the relations $[uvv_q,v^i(uvv_q)v^{-i}]=1$ for $i\in \{1,2,\ldots, k-1\}$; note also that $(uvv_q)^p=1$.

Finally, linear dependence of the $k+1$ vectors $\{\mathbf{e}_1,\ldots,\mathbf{e}_k,\mathbf{a}\}$ implies that $(\mathbf{0},I) = (\mathbf{a},I)\prod_{i=1}^k\ (\mathbf{e}_i,I)^{a_{i-1}}$. By \eqref{eq:ind} and \eqref{eq:a-k} this last equation translates to the relation
\begin{equation}\label{eq:long}
v^{k}(uvv_q)v^{-k}\prod_{i=0}^{k-1} \left(v^i(uvv_q)v^{-i}\right)^{a_{i}} = 1\ .
\end{equation}

For our primitive polynomial $\mu=\mu(\lambda)=\lambda^k+a_{k-1}\lambda^{k-1}+\ldots + a_1\lambda+a_0$, let $W(u,v;\mu)$ denote the left-hand side of \eqref{eq:long}; that is, we let
\begin{equation}\label{eq:W}
W(u,v;\mu) = v^{k}(uvv_q)v^{-k}\prod_{i=0}^{k-1} \left(v^i(uvv_q)v^{-i}\right)^{a_{i}}
\end{equation}
Collecting all the information gathered so far, we have identified the following relators in a presentation of the group $H_q=\langle u,v\rangle \cong {\rm AGL}(1,q)$ for an arbitrary prime power $q=p^k$, with $\mu$ and $W$ as above, and with $\ell=q$ if $q\in \{2,3\}$, $\ell=q-1$ if $q$ is even ($q\ne 2$) or $q\equiv 1$ mod $4$, and $\ell=(q-1)/2$ if $q\equiv 3$ mod $4$ but $q\ne 3$:
\begin{align}
\ & \mbox{$4$ relators determining orders of elements:}\ \ v^{q-1},\ (uv)^2,\ u^\ell,\ (uvv_q)^p, \label{eq:rel1} \\
\ & \mbox{$k-1$ commutativity relators:}\ \  [uvv_q,v^i(uvv_q)v^{-i}], \ 1\le i\le k-1, \ \ \mbox{and} \label{eq:rel2} \\  \ & \mbox{the relator implied by linear dependence:}\ \ W(u,v;\mu)\ . \label{eq:rel3}
\end{align}

Recall that $v_q=1$ for even $q$, and in this case the relator $(uvv_q)^p$ is redundant. There may be other redundancies for general $q$ but we will not investigate possibilities of further simplifications, except when $k=1$. In this case, there are no commutativity relators and the presentation found so far for $q=p$ and any primitive element $\xi \in F_p$, regarded as an integer mod $p$, reduces to
\begin{equation}\label{eq:k=1}
\langle u,v\ |\ u^\ell,v^{p-1},(uv)^2,(uv^{(p+1)/2})^p,v^{(p+1)/2}u(uv^{(p+1)/2})^\xi \rangle\ .
\end{equation}

We now show that the presentation given by \eqref{eq:rel1} -- \eqref{eq:rel3} determines a group isomorphic to ${\rm AGL}(1,q)$.

\begin{theorem}\label{thm:pres}
Let $p$ be an arbitrary prime, let $k$ be a positive integer, and let $\mu=\mu(\lambda)=\lambda^k +a_{k-1} \lambda^{k-1}+ \ldots + a_1\lambda +a_0$ be a monic primitive polynomial for a finite field of order $q=p^k$ with coefficients in the corresponding prime field. Further, let {\rm Rels} denote the collection of $k+4$ relators comprising the four relators from \eqref{eq:rel1}, the $k-1$ relators from \eqref{eq:rel2} and the single relator from \eqref{eq:rel3} that depends on $\mu$. Then, the group determined by the presentation $\langle u,v\ |\ {\rm Rels}\rangle$ is isomorphic to ${\rm AGL}(1,q)$.
\end{theorem}

{\bf Proof}. Let $H= H_q = \langle u,v\ |\ {\rm Rels}\rangle$. The derivation of the relators \eqref{eq:rel1} -- \eqref{eq:rel3} implies that ${\rm AGL}(1,q)$ is definitely a {\sl quotient group} of $H$; in particular, the subgroup generated by $v$ is isomorphic to $C_{q-1}$. Observe also that there is an epimorphism $\vartheta:\ H\to \langle v\rangle$ given by $\vartheta(u)=v_qv^{-1}$ and $\vartheta(v)=v$, or, equivalently, $\vartheta(uv)=v_q$ and $\vartheta(v)=v$; this follows from the fact that $\vartheta$ maps all the relators from the presentation of $H$ either to the identity or to the relator $v^{q-1}$ defining $C_{q-1}$. Let $L$ be the kernel of $\vartheta$; it follows that $H/L\cong C_{q-1}$ and $|H/L|=q-1$. Since ${\rm AGL}(1,q)$ is a quotient of $H$, the order of $L$ must be an integer multiple of $q$, and so our task will be accomplished if we show that $|L|=q$.
\smallskip

To do so, we will use the Reidemeister-Schreier method to determine a presentation for $L$. The method is well known and here we just adapt its brief description from \cite{Miy} to our situation, which simplifies because here $L$ is a normal subgroup of $H$.
\smallskip

We begin with setting up a Cayley digraph $D={\rm Cay}(H/L)$ for the group $H/L$ and the generating set $\{Lu,Lv\}$, which we may for simplicity identify with $\{u,v\}$; arcs of $D$ will be labelled by the corresponding generators and we will briefly refer to them as $u$-{\em arcs} and $v$-{\em arcs}. The $q-1$ vertices of $D$, each of in- and out-degree $2$, represent the $q-1$ right cosets $Lv^i$, $0\le i\le q-2$; instead of $Lv^i$ we will simply write $i$ in what follows. Observe that $D$ contains a directed Hamilton cycle formed by the $v$-arcs $(Lv^i, Lv^{i+1})$, or, in a simplified notation, $(i,i+1)_v$ for $i\in \{0,1,\ldots,q-2\}$, with $(i+1,i)_v$ being their reverses. As regards $u$-arcs, from $\vartheta(uvv_q)=1$ it follows that $uvv_q\in L$, and since $L\triangleleft H$ one also has $v^i(uvv_q)v^{-i}\in L$, $1\le i\le q-2$. The last containment relation is equivalent to $(Lv^i)u = Lv^{i-1}v_q$, which means that from every vertex $Lv^i$ there emanates a $u$-arc terminating at the vertex $Lv^{i-1}v_q$ which reduces to $Lv^{i+(q-3)/2}$ for $q$ odd (with exponents to be read mod $q-1$) and to $Lv^{i-1}$ for $q$ even. Letting $\overline{q}=(q-1)/2$ for $q$ odd and $\overline{q}=0$ for $q$ even, such $u$-arcs will be simply denoted $(i,i-1+\overline{q})_u$ for both parities of $q$. This completes the description of all arcs in our digraph $D$.
\smallskip

Observe that if $q$ is even, with $v_q=1$, this means that in the digraph $D$ for every $v$-arc $(i-1,i)_v$ one also has a $u$-arc $(i,i-1)_u$ in the direction reverse to $v$ (when $q=2$, the digraph has a single vertex with two loops). If $q=3$, the digraph has just $2$ vertices joined by a pair of mutually oppositely directed $v$-arcs, with a $u$-loop at each vertex. If $q=5$, then with every $v$-arc $(i,i+1)_v$ there is a parallel $u$-arc $(i,i+1)_u$. In all the remaining cases $D$ is a simple digraph, without loops and parallel arcs, containing another Hamilton cycle formed by $u$-arcs if $q\equiv 1$ mod $4$, and a pair of cycles of length $\overline{q}$ consisting of $u$-arcs if $q\equiv 3$ mod $4$.
\smallskip

The next step is choosing a (directed) spanning tree $T$ of $D$, which, in our case, will be formed by the Hamilton path consisting of the $q-2$ $v$-arcs $(i,i+1)_v$, $0\le i\le q-3$ (the tree will have a single vertex in the trivial case when $q=2$). With this choice of $T$, the $q$ cotree arcs comprise the $q-1$ $u$-arcs described above and the single $v$-arc $(q-2,1)_v$. For $0\le i\le q-2$, each $u$-arc $e_i$ from $Lv^i$ to $Lv^{i-1}v_q$ is contained in a unique closed walk in $T+e_i$, emanating from and terminating at the vertex $L$ (identified with $0$). For odd $q$ this walk has the form $P_iC_iP_i^{-1}$, where the closed walk $C_i$ (consisting just of a loop if $q=3$) and the open walk $P_i$ (with $P_0$ reducing to the vertex $0$) consist of arcs as follows:
\[ C_i = (i,i-1+\overline{q})_u\,(i-1+\overline{q},i-2+\overline{q})_v\,...\,(i+1,i)_v\ \ {\rm and} \ \  P_i = (0,1)_v\,...\, (i-1,i)_v\ .  \]
If $q$ is even, then the unique closed walk has the form $P_iC_iP_i^{-1}$, with $P_i$ as above and $C_i = (i,i-1)_u \,(i-1,i)_v$.
\smallskip

By Reidemeister-Schreier, for $i\in \{0,1,\ldots,q-2\}$, products of labels of arcs as they are encountered (in either direction) in the course of traversing the closed walk $P_iC_iP_i^{-1}$ constitute the so-called Schreier generators for $L$. In our case, if $q$ is odd, traversing the walk $P_iC_iP_i^{-1}$ gives the product $v^iu(v^{-1})^{ \overline{q}-1} v^{-i} = v^i(uvv_q)v^{-i}$, and if $q$ is even (and $v_q=1$), then the traversal yields the product $v^i(uv)v^{-i}$. It follows that for both parities of $q$ the $q-1$ $u$-arcs induce $q-1$ Schreier generators of the form $g_i=v^i(uv)v_qv^{-i}$ for $i\in \{0,1,\ldots,q-2\}$; note that they are mutually conjugate. The only $v$-arc not in $T$ is the arc $(v^{q-2},v)$ and because of $v^{q-1}=1$ it induces the identity element of $L$ as the last (and obviously redundant) Schreier generator. The important fact implied by Reidemeister-Schreier method is that $\{g_i\ |\ 0\le i\le q-2\}$ is a {\em generating set} for the group $L$.
\smallskip

It remains to determine new relators for a presentation of $L$. Normality of $L$ enables us to reduce the Reidemeister-Schreier process to considering only those closed walks representing relators in \eqref{eq:rel1} -- \eqref{eq:rel3} in the Cayley digraph $D$ which start (and end) at the vertex $L$, with the goal of rewriting each original relator in terms of Schreier generators. But since we know from the beginning that the order of $L$ is a multiple of $q$, it will suffice to apply the process to a suitable subset of original relators, namely, to the last relator of \eqref{eq:rel1} and to all of \eqref{eq:rel2} and \eqref{eq:rel3}. These can be rewritten in terms of our Schreier generators in a straightforward way by inspection, which results in the following collection of $k+1$ new relators: a single relator $g_0^p$ arising from the last relator in of \eqref{eq:rel1}, $k-1$ relators $[g_0,g_i]$ for $i\in \{1,2,\ldots, k-1\}$ arising from \eqref{eq:rel2}, and a single relator $g_k\prod_{i=0}^{k-1}(g_{i})^{a_{i}}$ coming from \eqref{eq:rel3} in which exponents are coefficients of the minimal polynomial $\mu$.
\smallskip

A repeated application of conjugation by $v$ to the relators $[g_0,g_i]$ for $i\in \{1,2,\ldots, k-1\}$ shows that  elements of the set $S_k=\{g_i\ |\ 0\le i\le k-1\}$ of the first $k$ mutually conjugate Schreier generators commute with each other, and each of them has order $p$ because of the relator $g_0^p$. The relator involving coefficients of $\mu$ enables one to write $g_k$ by means of generators in $S_k$, and by induction (again using repeated conjugation by $v$) it follows that each of the $q-1$ non-identity Schreier generators can be expressed in terms of elements of $S_k$. This shows that $L$ is Abelian and generated by $S_k$; as a consequence, $L$ has exponent $p$, and so $|L|\le p^k=q$. Since we have established that $|L|\ge q$ at the beginning, it follows that $|L|=q$ and $L\cong C_p^k$. (Note that here the set of all Schreier generators of $L$ found by the Reidemeister-Schreier method coincides with the entire group $L$.) \hfill $\Box$

\section{Automorphism groups of orientably-regular \\ embeddings of complete multigraphs}\label{sec:multi-compl}

Let $q=p^k$ be a prime power and let ${\cal M}(q,t)$ be an orientably-regular map with underlying graph $K_q^{(t)}$. In section \ref{sec:multi} we saw that for $q>2$ the group $G_{q,t}= \langle x,y\rangle \cong \Aut^+{\cal M}(q,t)$ is a central extension of the group $H_q\cong {\rm AGL}(1,q)$ by the cyclic group $N=\langle y^{q-1}\rangle$ of order $t$. We begin by a more general observation about such central extensions.
\smallskip

\begin{proposition}\label{prop:Schur}
For $q=p^k$ and $s\ge 1$ let $G$ be a group containing an element $y$ of order $(q-1)s$ such that $N = \langle y^{q-1}\rangle$ is a central subgroup of $G$, with $G/N\cong {\rm AGL}(1,q)$. If $q=4$, the commutator subgroup $G'$ of $G$ satisfies $|G' \cap N|\le 2$; if $q\ne 4$, then $G'\cong C_p^k$ and $G' \cap \langle y\rangle = 1$.
\end{proposition}

{\bf Proof.} By a result of \cite{Sch} on central and commutator subgroups of a group, the intersection $N\cap G'$ embeds in the Schur multiplier $M(G/N)\cong M({\rm AGL}(1,q))$. Luckily, the latter is known to be trivial for every $q\ne 4$, and isomorphic to $C_2$ if $q=4$, cf. \cite{Holt}. It follows that for $q\ne 4$ one has $|G'\cap N|=1$, implying that the normal subgroup $G'N$ of $G$ has order $|G'|s$ in this case, while $|G' \cap N|\le 2$ for $q=4$.
\smallskip

This enables us to specify $G'$ for $q\ne 4$. The commutator subgroup $(G/N)'$ of $G/N\cong {\rm AGL}(1,q)$ is known to be isomorphic to $C_p^k$, of order $q=p^k$. Combining the well-known formula $(G/N)'\cong G'N/N$ with the conclusion $|G'N| = |G'|s$ for $q\ne 4$ gives $q=|(G/N)'|=|G'N|/|N|=|G'|s/s=|G'|$, so that $|G'|=q$ if $q\ne 4$. Also, as $N\cap G'$ is trivial, from $(G/N)'\cong G'N/N$ it follows that $(G/N)'\cong G'$ and hence $G'\cong C_p^k$ if $q\ne 4$.
\smallskip

To finish the proof for $q\ne 4$, we extend the finding that $G'$ intersects trivially with $N$ to showing that $G' \cap \langle y\rangle = 1$; we may clearly also assume that $q\ne 2$. Suppose that $y^j\in G'$ for some positive integer $j < (q-1)s$ which is not divisible by $q-1$. Then, the order of $y^j$, which is equal to $(q-1)s/\gcd \{(q-1)s,j\}$, must divide $q=|G'|$. But since $q-1$ and $q$ are coprime, the denominator $\gcd\{(q-1)s,j\}$ must be divisible by $q-1$ and $j$ would have to be a multiple of $q-1$, a contradiction with triviality of $G'\cap N$. \hfill $\Box$
\bigskip

We will need Proposition \ref{prop:Schur} as a stand-alone statement later. As the next step we prove a more specific version of this result, related to the maps ${\cal M}(q,t)$ and the associated groups $G_{q,t}=\Aut^+{\cal M}(q,t)$. Here we will refer to the presentation of the quotient $H_q=\langle u,v\rangle$ of $G_{q,t}$ defined by means of the relators listed in \eqref{eq:rel1} -- \eqref{eq:rel3}.

\begin{proposition}\label{prop:Schur+}
Let $G = G_{q,t}=\langle x,y\rangle$ for $q = p^k$, with a central subgroup $N = \langle y^{q-1}\rangle \cong C_t$ such that $G/N\cong H_q=\langle u,v\rangle$, where $u=Nx$ and $v=Ny$, with presentation defined by {\rm \eqref{eq:rel1} -- \eqref{eq:rel3}} using a primitive polynomial $\mu$ over $F_q$. If $q \ne 4$, then $G'\cong C_p^k$ and $G' \cap \langle y\rangle = 1$. In particular, $G$ is then isomorphic to a split extension of the elementary abelian group $C_p^k$  by a cyclic group of order $(q-1)t$, in which the action of $\langle y\rangle$ on $C_p^k$ is given by $y^j\mapsto A^j$, where $A$ is the associated matrix of the primitive polynomial $\mu$ from {\rm \eqref{eq:m}}.
\end{proposition}

{\bf Proof.}
The first statement follows from Proposition \ref{prop:Schur} for $s=t$. This also implies that $G$ is a split extension of its subgroup $G'\cong C_p^k$ by the cyclic subgroup $\langle y\rangle < G$, of order $(q-1)t$. It remains to specify the conjugation action of $\langle y\rangle$ on $G'$ that determines the extension, that is, the corresponding group homomorphism $\theta:\ \langle y\rangle \to {\rm Aut}(C_p^k)\cong {\rm GL}(k,p)$. Since $N=\langle y^{q-1}\rangle$ is central in $G$, it follows that the value of $\theta$ on the coset $\bar{y}=Ny$ is constant. At the same time, the action of $\langle y\rangle$ on $G'$ in the group $G$ projects onto the action of $\langle \bar{y}\rangle$ on $G/N$ (which is still isomorphic to $G'$) in the factor group $G/N \cong H_q \cong C_p^k \rtimes \langle \bar{y}\rangle$. But the latter action is generated by $\bar{y}  \mapsto A$ for the matrix $A$ associated with $\mu$, which allows one to conclude that $\theta(y^j) =A^j$ for every $j$ mod $(q-1)t$. \hfill $\Box$
\bigskip

While the extension $C_p^k\rtimes C_{(q-1)t}$ for $q\ne 4$ is uniquely determined {\em as an abstract group} (for a given matrix $A$), this does not necessarily mean that it uniquely determines a map ${\cal M}(q,t)$ with ${\rm Aut}^+ {\cal M}(q,t)=\langle x,y\rangle$. Theorem \ref{thm:Gpt} and Proposition \ref{prop:q=4} will show that there may exist several {\em distinct} choices of elements $x$ with the same coset $\bar{x}=Nx$ but giving non-isomorphic maps defined by distinct generating pairs $(x,y)$ for $y$ as in Proposition \ref{prop:Schur+}.
\smallskip

By the theory outlined in the introduction, the map ${\cal M}(q,t)$ is a $t$-fold branched cover of an orientably-regular embedding ${\cal M}_q$ of a {\sl simple} complete graph of order $q$, and if $q>2$, isomorphism classes of the maps ${\cal M}_q$ are in one-to-one correspondence with primitive polynomials over $F_q$. Choose such a polynomial $\mu=\mu(\lambda)=\lambda^k +a_{k-1}\lambda^{k-1}+\ldots + a_1\lambda +a_0$ and, as before, let $ W(u,v;\mu)$ be given by \eqref{eq:W}. Also, we keep $\ell$ equal to $q-1$ if $q$ is even or $q\equiv 1$ mod $4$, and equal to $(q-1)/2$ if $q\equiv 3$ mod $4$, save the exceptional case $\ell=3$ for $q=3$. Theorem \ref{thm:pres} combined with the `lifting' method from \eqref{eq:lift} imply that for $q>2$ the group $G_{q,t}=\Aut^+{\cal M}(q,t)$ admits a presentation of the form
\begin{equation}\label{eq:Gqt}
G_{q,t}=\langle x,y\ |\ y^{(q{-}1)t},\,(xy)^2,\,[x,y^{q{-}1}],\, {\rm Rel}_1,\,{\rm Rel}_2,\,{\rm Rel}_3\, \rangle
\end{equation}
where the collections of relators ${\rm Rel}_1$, ${\rm Rel}_2$ and ${\rm Rel}_3$ are obtained, respectively, from the relators in \eqref{eq:rel1}, \eqref{eq:rel2} and \eqref{eq:rel3} as follows. First, we skip the first two relators of \eqref{eq:rel1} because their new versions appear in \eqref{eq:Gqt}. Then, we substitute $x$ for $u$, $y$ for $v$, and $y_q=y^{(q-1)/2}$ for the value of $v_q$ if $q$ is odd (with $y_q=v_q=1$ for $q$ even) in each of the remaining relators of \eqref{eq:rel1} -- \eqref{eq:rel3}. Finally, to each of the resulting words (in $x$ and $y$) arising from such substitutions we append a multiplicative factor consisting of some power of $y^{1-q}$. In our notation, the multiplicative factors will have the form $y^{(1-q)\omega}$ for $\omega\in \{a,b,c,e_1,\ldots,e_{k-1}\}$ mod $t$, and hence the collection of relators in the presentation \eqref{eq:Gqt} will be as follows:
\begin{align}
\ & \mbox{${\rm Rel}_1$:}\ \  x^\ell\cdot y^{(1{-}q)a},\  (xyy_q)^p\cdot y^{(1{-}q)b}; \label{eq:rel1t} \\
\ & \mbox{${\rm Rel}_2$:}\ \  [xyy_q,y^ixyy_qy^{-i}] \cdot y^{(1{-}q)e_i}, \ 1\le i\le k-1;  \label{eq:rel2t} \\
\ & \mbox{${\rm Rel}_3$:}\ \ W(x,y;\mu)\cdot y^{(1{-}q)c}\ . \label{eq:rel3t}
\end{align}
Recall that $y_q=1$ for even $q$; note also that the value of $f$ appearing in \eqref{eq:lift} is here equal to $1$ since our underlying graphs for $q>2$ are not bipartite.
\smallskip

We will determine the parameters $a$, $b$, $c$ and $e_i$ ($1\le i\le k-1$) mod $t$, under the assumptions that $q>4$ and $|G_{q,t}|=q(q-1)t$, which is the expected order of the group $\Aut^+ {\cal M}(q,t)$.

\begin{theorem}\label{thm:Gpt}
For a prime power $q=p^k>4$ let $G = \langle x,y\rangle$ be a group with a presentation given by {\rm \eqref{eq:Gqt} -- \eqref{eq:rel3t}}, and for odd $p$ let $\alpha=\alpha(\mu)$ be the unique even integer in the set $\{1,\ldots, p-1\}\cup\{p+1,\ldots,2p-1\}$, equal {\rm mod} $p$ to the sum of the coefficients of the primitive polynomial $\mu$. If the order of $G$ is $q(q-1)t$, then the following holds:
\medskip

\noindent {\rm (I)} If $q$ is odd, then $t$ is also odd, and {\rm mod} $t$ one has $e_i=0$ {\rm ($1\le i\le k-1$)}, $a=-1$ for $q\equiv 1$ {\rm mod} $4$ while $a=(t-1)/2$ for $q\equiv 3$ {\rm mod} $4$, $b=(p+t)/2$, and $c=\alpha/2$.
\medskip

\noindent {\rm (II)} If $q$ is even, then  $b=e_i=0$ {\rm ($1\le i\le k-1$)} for every $t\ge 1$; moreover, if $t$ is odd, then $(a,c)=(-1,0)$, and if $t$ is even, then $(a,c)\in \{(-1,0),(t/2-1,t/2)\}$, all {\rm mod} $t$.
\medskip

\noindent Conversely, if $G$ is a group with presentation {\rm \eqref{eq:Gqt} -- \eqref{eq:rel3t}} in which $q>4$ and the parameters $t$, $a$, $b$, $c$ and $e_i$ satisfy {\rm (I)} for $q$ odd and {\rm (II)} for $q$ even, then $G$ has order $q(q-1)t$.
\end{theorem}

{\bf Proof.} Let $G$ be a group presented as in \eqref{eq:Gqt} and assume that its order is $q(q-1)t$. Clearly, $N=\langle y^{q-1}\rangle$ is a central subgroup of $G$; by Theorem \ref{thm:pres} and the assumption on the order of $G$ one has $|N|=t$, and the order of $y$ is $(q-1)t$. Further, we know that the quotient $G/N$ is isomorphic to $H_q\cong {\rm AGL}(1,q)$. For $q > 4$, Proposition \ref{prop:Schur+} implies that $|G'|=q$ and $G'\cap \langle y\rangle$ is trivial; hence the normal subgroup $G'N$ of $G$ has order $|G'|t$.
\smallskip

Thus, the abelianisation $G/G'$ is isomorphic to $\langle y\rangle\cong C_{(q-1)t}$. In what follows we let $\bar{x}=G'x$ and identify $y$ with $G'y$ (which we may do as no non-trivial power of $y$ is in $G'$); also, we let $i$ extend over the set $\{1,\ldots, q-1\}$. Our first version of the induced presentation of $G/G'$, obtained only by absorbing commutators in \eqref{eq:rel2t}, has the form
\begin{equation}\label{eq:G/G'}
\langle \bar{x},y\ |\ y^{(q{-}1)t},\,(\bar{x}y)^2,\, (\bar{x})^\ell y^{(1{-}q)a},\, (\bar{x}yy_q)^py^{(1{-}q)b},\, \ldots,\, y^{(1-q)e_i}, \,\ldots,\, \, W(\bar{x},y;\mu)y^{(1-q)c}\,    \rangle
\end{equation}

The facts that the order of $y$ in \eqref{eq:G/G'} is {\em equal to} $(q-1)t$ (by our assumption that $|G_{q,t}| =q(q-1)t$) together with $[\bar{x},y]=1$ enable us to determine the parameters $a,b,c$ and $e_i$ ($1\le i\le k-1$) from \eqref{eq:G/G'} quite conveniently. For example, because of the order of $y$, presence of the relators $y^{(1-q)e_i}$ in \eqref{eq:G/G'} implies that $(q-1)t\mid (q-1)e_i$ and hence $e_i\equiv 0$ mod $t$ for every $i\in \{1,\ldots, k-1\}$. The remaining parameters, in the order $b$, $a$ and $c$, will be calculated in the next three steps.
\smallskip

{\sl Step 1: the value of $b$}. Observe that the relator $(\bar{x}yy_q)^py^{(1{-}q)b}$ for even $q$ collapses to $y^{(1{-}q)b}$, implying that $b\equiv 0$ mod $t$, with no restriction on the parity of $t$. For odd $q$, with the help of $(\bar{x}y)^p=\bar{x}y$ in $G/G'$ the same relator reduces to $(\bar{x}y)\cdot y^{p(q-1)/2}\cdot y^{(1-q)b}$ and hence gives $\bar{x}y=y^{(2b-p)(q-1)/2}$. Since $G'$ has odd order for odd $q$, the element $\bar{x}y$ has order $2$ in $G/G'$ and so $y^{(2b-p)(q-1)/2}$ must be an involution in the group $G/G'\cong \langle y\rangle$. As the latter contains a single involution $y^{(q-1)t/2}$, it follows that $(q-1)t/2$ divides $(2b-p)(q-1)/2$ and so $t\mid (2b-p)$. In particular, if $q$ is odd, then so is $t$, and one may take $b=(p+t)/2$.
\smallskip

{\sl Step 2: the value of $a$}. To calculate $a$ from the relator $(\bar{x})^\ell y^{(1{-}q)a}$, we begin with odd $q$. If $q\equiv 1$ mod $4$, then $\ell=q-1$ is even and so $(\bar{x}y)^{q-1}=1$ in $G/G'$, which reduces the relator to $y^{(1{-}q)(a+1)}$, and by the order of $y$ we conclude that $a\equiv -1$ mod $t$. If $q\equiv 3$ mod $4$, then $\ell=(q-1)/2$ is odd and the relator is equivalent to $(\bar{x})^{(q-1)/2}y^{(1-q)a} = (\bar{x}y)^{(q-1)/2}\cdot y^{(1-q)a + (1-q)/2} = \bar{x}y\cdot y^{(1-q)a + (1-q)/2}$, and using the fact that $\bar{x}y$ is an involution we obtain that $(q-1)t/2$ divides $(2a+1)(q-1)/2$ and so $2a+1\equiv 0$ mod $t$ (and one may take $a=(t-1)/2$, knowing that $t$ must be odd here). If $q$ is even, then $\ell = q-1$ is odd, and $(\bar{x}y)^{q-1} = \bar{x}y$ in $G/G'$. Now, $\bar{x}y$ has order $1$ or $2$ in $G/G'$, depending on whether $xy\in G'$ or not (and both cases occur, as we shall see later). Accordingly, the relator we started with reduces to $y^{(1{-}q)(a+1)}$ if $xy\in G'$, and to $\bar{x}y\cdot y^{(1{-}q)(a+1)}$ if $xy\notin G'$. This implies that $t\mid (a+1)$ in the first case, and $t\mid 2(a+1)$ but $t\nmid (a+1)$ in the second case (for $t$ even), by calculations similar to the previous ones; in particular, one may take $a=-1$ and $a=t/2-1$, respectively.
\smallskip

{\sl Step 3: the value of $c$}. The relator $W(\bar{x},y;\mu) = y^k(\bar{x}yy_q)y^{-k}\prod_i y^i(\bar{x}yy_q)^{a_i }y^{-i} \cdot y^{(1-q)c}$ will enable us to determine $c$ by a successive abelianisation. First, canceling conjugates by powers of $y$ gives  $(\bar{x}yy_q)^{a(\mu)}y^{(1-q)c}$, where $\alpha = \alpha(\mu)$. Clearly, $\alpha$ is not congruent to $0$ mod $p$, because of irreducibility of $\mu$. If $q$ is even, that is, if $p=2$, then the last relator collapses to $\bar{x}y\cdot y^{(1-q)c}$, leading to the conclusion that $t\mid c$ if $xy\in G'$ and $t\mid 2c$ but $t\nmid c$ if $xy\notin G'$ (in which case $t$ must be even). This means that if $xy\in G'$, one may take $c=0$ for every $t\ge 1$, and if $xy\notin G'$, the element $xy$ has order $2$ in $G/G'$, forcing $t$ to be even, and then one may take $c=t/2$. If $q$ is odd, then $\bar{x}y$ has order $2$ in $G/G'$, because $|G'|=q$ is odd. As $\alpha$ mod $t$ is the unique {\em even} integer between $2$ and $2t-2$, the relator $(\bar{x}yy_q)^{\alpha} y^{(1-q)c}$ further reduces to $y^{\alpha(q-1)/2}y^{(1-q)c}= y^{(\alpha-2c) (q-1)/2}$, giving $t\mid (\alpha/2 - c)$, and one may take $c=\alpha/2$. (The reader may have noticed that different choices of integers $\alpha=\alpha(\mu)$ mod $t$ for odd $q$ would give different values of $c=c(\alpha)$, but in conjunction with the relator $(\bar{x}yy_q)^py^{(1-q)b}$ it may be checked that any pair of such choices determine equivalent presentations.)
\smallskip

For the converse, suppose that $G$ is given by the presentation {\rm \eqref{eq:Gqt} -- \eqref{eq:rel3t}} in which the parameters $t$, $a$, $b$, $c$ and $e_i$ satisfying {\rm (I)} for $q$ odd and {\rm (II)} for $q$ even. Let us now substitute these values in the corresponding relators of \eqref{eq:G/G'} that define the abelianisation $G/G'$, referring in the process to our calculations in Steps 1 -- 3. To begin with,  note that that all the relators in \eqref{eq:G/G'} containing $e_i$ simply vanish by letting $e_i=0$, $1\le i\le k-1$.
\smallskip

By Step 1, the relator containing $b$ vanishes for $b=0$ if $q$ is even, and gives $\bar{x}y=y^{t(q-1)/2}$ for $b=(p+t)/2$ if $q$ is odd, which implies that $\bar{x}=y^{t(q-1)/2-1} \in \langle y\rangle$. By Step 2, the relator involving the parameter $a$ vanishes for $q\equiv 1$ mod $4$ (when $a=-1$) and gives $\bar{x}y= y^{t(q-1)/2}$ for $q\equiv 3$ mod $4$ (when $a=(t-1)/2$), so that $\bar{x}=y^{t(1-q)/2-1}$ again. For $q$ even the analysis depended on whether $xy\in G'$ or not; if $xy\in G'$, then $\bar{x}y=1$ in $G/G'$, so that $\bar{x}=y^{-1}$ (for $a=-1$ and every $t\ge 1$), and if $xy\notin G'$, then one arrives at $\bar{x}=y^{(q-1)(t/2)-1}$ for $a=t/2-1$. Finally, by Step 3 the relator containing $c$ either vanishes or implies that $\bar{x}=y^{t(q-1)/2-1}$, by reasoning analogous to those related to Step 1 and 2. In all these cases the formulae for $\bar{x}$ show not only that $\bar{x}\in \langle y\rangle$ but also that $(\bar{x}y)^2=(\bar{x})^2y^2=1$.
\smallskip

It follows that for $a$, $b$, $c$ and $e_i$ given by (I) and (II), the presentation \eqref{eq:G/G'} of the group $G/G'$ collapses to $\langle y\ |\ y^{(q-1)t}\rangle$.
\smallskip

But observe that by the presentation \eqref{eq:Gqt} -- \eqref{eq:rel3t} of the original group $G$, the subgroup $N=\langle y^{q-1}\rangle$ is central in $G$, and one may check against the list of relators in \eqref{eq:rel1} -- \eqref{eq:rel3} that the induced presentation of $G/N$ determines a group isomorphic to ${\rm AGL}(1,q)$. Applying Proposition \ref{prop:Schur} to this situation for $q > 4$ gives $G'\cap \langle y\rangle =1$ and $|G'|=q$. Taken together with the derived induced presentation of $G/G'$ one concludes that the order of $y$ is {\em equal} to $(q-1)t$ and hence $|G|=q(q-1)t$, as claimed. This completes the proof.   \hfill $\Box$
\bigskip

Thus, for an arbitrary prime power $q=p^k> 4$ and an arbitrary integer $t\ge 1$, a complete graph  $K_q^{(t)}$ of order $q$ and edge-multiplicity $t$ admits an orientably-regular embedding ${\cal M}$ if and only if $t$ and the parameters $a$, $b$, $c$ and $e_i$ ($0\le i\le k-1$) appearing in {\rm \eqref{eq:Gqt} -- \eqref{eq:rel3t}} are as stated in parts (I) and (II) of Theorem \ref{thm:Gpt}, in which case the group ${\rm Aut}^+({\cal M})$, of order $q(q-1)t$, admits a presentation given by \eqref{eq:Gqt} -- \eqref{eq:rel3t}. For a given $q>4$ and a monic primitive polynomial $\mu$ over $F_q$, this gives two maps if both $q$ and $t$ are even, which we will denote ${\cal M}_0(q,t;\mu)$ and ${\cal M}_1(q,t;\mu)$, where the subscripts correspond, respectively, to the values of $c=0$ and $c=t/2$ from part (II) of Theorem \ref{thm:Gpt}, and a unique map for every odd $t\ge 1$ and an arbitrary $q> 4$, which we will denote ${\cal M}_2(q,t;\mu)$. It remains to check for potential isomorphisms between these maps.
\smallskip

\begin{theorem}\label{thm:iso}
Let $S$ be the set of ordered triples $(q,t;\mu)$, where $q>4$ is a prime power, $t$ is a positive integer {\rm (}odd if $q$ is odd{\rm)}, and $\mu$ is a monic primitive polynomial over $F_q$. Among all maps ${\cal M}_i (q,t;\mu)$ for $(q,t;\mu)\in S$ and $i\in \{0,1,2\}$, no two are isomorphic to each other. In particular, the number of pairwise non-isomorphic orientably-regular embeddings of $K_q^{(t)}$ for $q=p^k>4$ is equal to $\varphi(q-1)/k$ if $t$ is odd, and $2\varphi(q-1)/k$ if both $q$ and $t$ are even.
\end{theorem}

{\bf Proof.} Let $(q,t;\mu)$ and $(q',t';\mu')$ be two triples from $S$, and let ${\cal M}={\cal M}_i(q,t;\mu)$ and ${\cal M}'={\cal M}_j(q',t';\mu')$ be maps as in the statement, for some $i,j\in \{0,1,2\}$. Suppose that the map ${\cal M}$ is isomorphic to the map ${\cal M}'$. Then, their automorphism groups have the same order, so that $q(q-1)t = q'(q'-1)t'$, and they have the same valency, that is, $(q-1)t=(q'-1)t'$; the two equalities imply that $q=q'$  and hence $t=t'$. It follows that the quotients of the two maps obtained by factoring out the central cyclic subgroup of order $t=t'$ are also isomorphic. But these quotients are orientably-regular embeddings of a (simple) complete graph $K_q$, and by \cite{JaJo} their isomorphism is equivalent to $\mu=\mu'$.
\smallskip

To finish the proof it remains to distinguish, for a fixed triple $(q,t,\mu)$ with even $q>4$, between the maps ${\cal M}_0(q,t;\mu)$ and ${\cal M}_1(q,t;\mu)$. For $i\in \{0,1\}$ let $G_i = {\rm Aut}^+{\cal M}_i(q,t;\mu) = \langle x_i,y_i\rangle$ be the corresponding groups, with presentations \eqref{eq:Gqt} -- \eqref{eq:rel3t} obtained by substituting $x_i$ and $y_i$ for $x$ and $y$. By \cite{JoSi}, the two maps are isomorphic if and only if the assignment $(x_0,y_0)\mapsto (x_1,y_1)$ extends to a group isomorphism $G_1\to G_2$. This, however, is excluded by the relator $x^{\ell}y^{(1-q)a}$ in \eqref{eq:rel3t} and the distinct values of $a=-1+it/2$ for $i\in \{0,1\}$.
\smallskip

The claim about the number of such embeddings follows immediately from Theorem \ref{thm:Gpt} and the known number of distinct monic primitive polynomials over $F_q$. \hfill $\Box$
\bigskip

Types $\{m,n\}$ of the maps from Theorem \ref{thm:iso} are given by $m={\rm ord}(x)$ and $n={\rm ord}(y)$ for $G=\langle x,y\rangle$ as in \eqref{eq:Gqt} -- \eqref{eq:rel3t}. We had $n=(q-1)t$ throughout, and $m$ will now be determined from the relation $x^\ell = y^{(q-1)a}$, implied by \eqref{eq:rel1t}. The key is to observe that, taking the quotient of $G/N\cong {\rm AGL}(1,q)$ by the central subgroup $N=\langle y^{q-1}\rangle$ gives ${\rm ord}(Nx)=\ell$, so that $\ell$ is the smallest positive integer such that $x^\ell\in \langle y^{q-1}\rangle$. It follows that ${\rm ord}(x)$ is equal to $\ell$ times the order of $y^{q-1}$ in $\langle y\rangle$, that is, ${\rm ord}(x) = \ell t/\gcd(a,t)$. Working out the values of $\gcd(a,t)$ for the values of $q$, $a$ and $t$ from Theorem \ref{thm:Gpt} gives the following table of types of orientably-regular embeddings of $K_q^{(t)}$ for $q>4$:
\begin{table}[hbt!]
	\centering
\begin{tabular}{c|c|c|c|c}
Prime power          &  Multiplicity   &  Value of $a$        & Map         &  Type \\ \hline
$q=2^k > 4$    &  every $t\ge 1$       &  $-1$       & \ ${\cal M}_0(q,t;\mu)$ \ & $\{(q-1)t,\, (q-1)t\}$ \\
$q=2^k > 4$    &  $t\equiv 0$ mod $4$ &  $t/2{-}1$  & \ ${\cal M}_1(q,t;\mu)$ \ & $\{(q-1)t,\, (q-1)t\}$ \\
$q=2^k > 4$    &  $t\equiv 2$ mod $4$  &  $t/2{-}1$  & \ ${\cal M}_1(q,t;\mu)$ \ & $\{(q-1)t/2,\, (q-1)t\}$ \\
$q\equiv 1$ mod $4$  &  odd $t\ge 1$   &  $-1$       & \ ${\cal M}_2(q,t;\mu)$ \ & $\{(q-1)t,\,(q-1)t\}$ \\
$q\equiv 3$ mod $4$  &  odd $t\ge 1$   &  $(t{-}1)/2$  & \ ${\cal M}_2(q,t;\mu)$ \ & $\{(q-1)t/2,\,(q-1)t\}$ \\
\end{tabular}
\caption{Types of $t$-inflations of complete orientably-regular maps for $q>4$.} \label{tab:q-gen}
\end{table}

The genus $g$ of an orientably-regular map ${\cal M}$ of type $\{m,n\}$ with $G={\rm Aut}^+{\cal M}$ is given by Euler's formula. Our maps have $|G|/m$ vertices, $|G|/2$ edges and $|G|/n$ faces, so that the formula reads $|G|/m-|G|/2+|G|/n = 2-2g$. For maps in Table \ref{tab:q-gen}, by evaluating the genus for $|G|=q(q-1)t$ one obtains  $g=1-q+q(q-1)t/4$ for type ${m,n}$ with $m=n=(q-1)t$, and $g= 1-3q/2 + q(q-1)t/4$ for type $\{m,n\}$ where $m=(q-1)t/2$ and $n=(q-1)t$.

\section{The three smallest prime powers}\label{sec:sing}

For the cases $q\in \{2,3,4\}$ avoided so far one could simply reproduce available results. Our aim, however, is to re-prove them by applying the method of `lifting presentations' and show how the associated irregularities can be overcome.
\smallskip

The case $q=2$ represents a singularity not only by its bipartite underlying graphs on $2$ vertices with $t$ parallel edges (a $t$-{\em dipole}) but also by face length $\ell=q=2$ for $K_2$ on a sphere. Orientably-regular embeddings of $t$-dipoles have been classified in \cite{Gar} and later in \cite{NeSk}, but the results follow easily also from our approach outlined in section \ref{sec:multi}. Indeed, by the general form \eqref{eq:lift} of the lift of a presentation, orientably-regular embeddings of $t$-dipoles are in a one-to-one correspondence with group presentations of the form
\begin{equation}\label{eq:q=2}
G_{2,t} =\langle\, x,y\ |\ x^m,\,y^t,\,(xy)^2,\,x^2y^{f+1}\rangle\ \ {\rm with}\ \  f^2\equiv 1\ {\rm mod\ } t
\end{equation}
where $m=2t/\gcd(f+1,t)$ is the face length, and hence in a one-to-one correspondence with solutions $f$ mod $t$ of the congruence $f^2\equiv 1$ mod $t$. The number of solutions of this congruence (and so the number of orientably-regular embeddings of a $t$-dipole) is known to be $2^{r+s}$, where $r$ is the number of distinct odd prime factors of $t$ and $s=2,\, 1$ or $0$ depending on whether $t$ is divisible by $8$, by $4$ but not by $8$, and not divisible by $4$, respectively; the genus of the embedding for such an $f$ mod $t$ evaluates to $(t-\gcd(f{+}1,t))/2$.
\smallskip

The value $q=3$ is exceptional because of the face length $\ell=q=3$ for a triangle on a sphere. Orientably-regular embeddings of $K_3^{(t)}$ have been classified in \cite{Gar}, and later also in \cite{H+} as part of a study of inflations of embedded cycles. The result easily follows by replacing $\ell$ with $3$ in Step 2 of the proof of Theorem \ref{thm:Gpt}, giving $a=(t-3)/2$ for $q=3$. It follows that orientably-regular embeddings of $K_3^{(t)}$ are in a one-to-one correspondence with presentations $\langle\, x,y\ |\ x^m=y^n=(xy)^2=[x,y^2]=x^3y^{3-t}=1 \,\rangle$ for $t$, $m$ and $n$ (and genus) as in Table \ref{J:tab:3vert}, with a unique embedding for every odd $t\ge 1$.
\begin{table}[hbt!]
	\centering
\begin{tabular}{c|c|c}
Multiplicity  & Type $\{m,n\}$ & Genus   \\ \hline
$t$ odd, $3\mid t$ & $\{t,2t\}$ &  $(3t-7)/2$  \\
$t$ odd, $3\nmid t$ & $\{3t,2t\}$ & $(3t-3)/2$ \\
\end{tabular}
\caption{Orientably-regular embeddings of $K_3^{(t)}$.} \label{J:tab:3vert}
\end{table}

Finally, let $q=4$, a case solved in \cite{Gar} as well; it stands out by allowing a non-trivial Schur multiplier in Proposition \ref{prop:Schur}. We begin by simplifying the presentation \eqref{eq:Gqt} -- \eqref{eq:rel3t} of the group $G=G_{q,t}$ for $q=4$, observing first that for $p=k=2$, $y_q=1$ and $\ell=q-1=3$ the relator in \eqref{eq:rel1t} containing $b$ is immaterial. As the only primitive polynomial for $q=4$ is $\lambda^2+\lambda+1$, one has $W(x,y;\mu)=y^2(xy)y^{-2}\cdot y(xy)y^{-1}\cdot (xy) = y^{2}x^3y$, and using $x^3=y^{3a}$ gives $W(x,y;\mu) y^{-3c}=y^{3(1+a-c)}$. Since the order of $y$ is the valency of the associated orientably-regular embedding of $K_4^{(t)}$, equal to $(q-1)t=3t$, we {\em assume} that ${\rm ord}(y)=3t$ in \eqref{eq:Gqt}, so that the relator $y^{3(1+a-c)}$ may be replaced by the condition $c=a+1$ mod $t$. Presentation of $G$ then simplifies as follows,  with $e_1$ replaced by $e$ for $k=2$ in \eqref{eq:rel2t}:
\begin{equation}\label{eq:q=4}
G=G_{4,t} = \langle x,y\ |\ x^3y^{-3a},\,y^{3t},\,(xy)^2,\,[x,y^3],[xy,yx]y^{-3e} \rangle
\end{equation}

We are now in position to clarify the situation for $q=4$, including an explanation of connections of this case with the Schur multiplier from Propositions \ref{prop:Schur} and \ref{prop:Schur+}.

\begin{proposition}\label{prop:q=4}
Let $G=\langle x,y\rangle$ be a group presented as in {\rm \eqref{eq:q=4}}, with $y$ of order $3t$, $t\ge 1$. If $|G|=12t$, then the parameters $a$ and $e$ {\rm mod $t$} in \eqref{eq:q=4} satisfy ${\rm (A)}$ or ${\rm (B)}$ below:
\medskip

\noindent{\rm (A)} $a=-1$ for arbitrary $t$, together with $a=-1+t/2$ for every even $t$; in both cases $e=0$, $G'= \langle [x^{-1},y],\,[y^{-1},x]\rangle \cong C_2\times C_2$ and $G'\cap \langle y\rangle = 1$;
\medskip

\noindent{\rm (B)} $(a,e)=(-1\pm t/4,t/2)$ for every $t\equiv 0$ {\rm mod} $4$, and $G'= \langle [x^{-1},y], \,[x,y],\,[y^{-1},x]\rangle$ is isomorphic to the quaternion group in both cases, with $G'\cap \langle y\rangle = \langle y^{3t/2}\rangle$.
\medskip

\noindent Conversely, if $G$ is a group with presentation {\rm \eqref{eq:q=4}} in which $t$, $a$ and $e$ satisfy $(A)$ or $(B)$, then $G$ has order $12t$. Moreover, the corresponding orientably-regular embeddings of $K_4$ with edge-multiplicity $t$ are mutually non-isomorphic, and their number is $1$, $2$, and $4$, depending on whether $t$ is odd, $t\equiv 2$ mod $4$, and $t\equiv 0$ mod $4$, respectively.
\end{proposition}

{\bf Proof.} Let $G$ be as in the statement. Abelianisation of the presentation \eqref{eq:q=4}, with $x_{\rm o}=G'x$ and $y_{\rm o}=G'y$, gives $G/G' = \langle x_{\rm o},\,y_{\rm o}\ |\ x_{\rm o}^3y_{\rm o}^{-3a},\,y_{\rm o}^{3t}, \,x_{\rm o}^2y_{\rm o}^2,\,y_{\rm o}^{-3e}\rangle $. Proposition \ref{prop:Schur} implies that the intersection $G' \cap \langle y^3\rangle$ has order at most $2$, and in case of equality one has $G' \cap \langle y^3\rangle = \langle y^{3t/2} \rangle$ for even $t$, as $y^{3t/2}$ is then the only involution in $\langle y\rangle \cong C_{3t}$. It follows that the order of $y_{\rm o}$ is divisible by $3t/2$. Comparing the values of $x_{\rm o}^6$ obtained by squaring the first relator and cubing the third relator of the presentation of $G/G'$ gives $y_{\rm o}^{6a} = y_{\rm o}^{-6}$, which implies that $3t/2$ divides $6(a+1)$ and hence $t\mid 4(a+1)$. This gives four values of $a$ mod $t$: $a=-1$ for any $t$, $a=-1+t/2$ for $t$ even, and $a=-1\pm t/4$ for $t\equiv 0$ mod $4$. Before beginning to consider cases, observe that centrality of both $x^3$ and $y^3$ from \eqref{eq:q=4} together with $(xy)^2=1$ imply that $(x^3y^3)^2=[xy,yx]$, independently of the value of $a$.
\smallskip

Let $a\in \{-1,-1+t/2\}$, that is, let $a=-1+\veps t/2$ for $\veps\in \{0,1\}$; the first relator of \eqref{eq:q=4} then gives $(x^3y^3)^2=1$. This means that also $[xy,yx]=1$, so that $y^{3e}=1$ in \eqref{eq:q=4}, and by ${\rm ord}(y)=3t$ one has $e=0$. Further, from $y^{3\veps t/2}=y^{3(a+1)}=x^3y^3= x\cdot[x^{-1},y]\cdot y$ one has $[x^{-1},y] = yxy^{3\veps t/2}$ and, similarly, $[y^{-1},x] = xyy^{3\veps t/2}$.
\smallskip

The relation $x^3y^3=y^{3\veps t/3}$ is equivalent to $y(yxy^{3\veps t/2})y^{-1} = (yx)(xy)$. The latter with $y(xyy^{3\veps t/2})y^{-1}=yxy^{3\veps t/2}$ implies that $G_0=\langle xyy^{3\veps t/2}, yxy^{3\veps t/2}\rangle = \langle [x^{-1},y],\,[y^{-1},x]\rangle$ is a normal subgroup of $G$, isomorphic to $C_2\times C_2$. Since $G=\langle x,y\rangle$, one has $G=G_0 \rtimes \langle y\rangle$, and it can be checked that $G_0 \cap \langle y\rangle=1$, so that $|G|=12t$. But the semidirect product structure of $G$ also implies that $G'=G_0$, with $G' \cap \langle y\rangle = 1$. This proves (A).
\smallskip

To deal with (B), let us substitute $a=-1\pm t/4$ in the first relator of the presentation \eqref{eq:q=4} of the original group $G$; this results in $x^3y^3=y^{\pm 3t/4}$ and $(x^3y^3)^2= y^{3t/2}$. But we saw earlier that $(x^3y^3)^2=[xy,yx]$, which now implies that $y^{3t/2}\in G'$. Since ${\rm ord}(y)=3t$ in $G$, the commutator $[xy,yx]$ is not trivial, and so comparison of the relations $(x^3y^3)^2=y^{3t/2}$ and $[xy,yx]=y^{3e}$ gives $e=t/2$ mod $t$. It may be checked that in this case the relator containing $e$ is a consequence of the remaining ones in \eqref{eq:q=4}.
\smallskip

Centrality of $x^3$ and $y^3$ with $(xy)^2=1$ also imply that $(x^3y^3)^2 =(x^2y^2)^2=[x^{-1},y]^2$, so that now one has $y^{3t/2}= (x^3y^3)^2=[x^{-1},y]^2$. It follows that the commutator $[x^{-1},y]$ has order $4$ in $G$, and one may show in a similar way that the commutators $[x,y]$ and $[y^{-1},x]$ have order $4$ as well. But the product of the last three commutators (in their order of appearance) is easily verified to be equal to the involution $y^{3t/2}$, implying that in this exceptional case $G'$ is isomorphic to the quaternion group. This establishes (B).
\smallskip

Existence of groups $G$ of order $12t$ with presentation \eqref{eq:q=4} and satisfying (A), (B) and the claim about their numbers can be shown as in the proofs of Theorems \ref{thm:Gpt} and \ref{thm:iso}. \hfill $\Box$
\bigskip

Types $\{m,n\}$ of the orientably-regular embeddings of the graphs $K_4^{(t)}$ for $q=4$ are determined by the same method that was used to set up Table \ref{tab:q-gen}. Here we just give a summary, including genera of the maps as a consequence of Euler's formula:

\begin{table}[hbt!]
	\centering
\begin{tabular}{c|c|c|c}
Multiplicity           & Value of $a$      &  Type          & Genus \\ \hline
any $t\ge 1$           & $-1$     & $\{3t,3t\}$    & $3t-3$ \\
$t\equiv 0$ mod $4$    & $t/2-1$  & $\{3t,3t\}$    & $3t-3$ \\
$t\equiv 0$ mod $8$    & $t/4-1$  & $\{3t,3t\}$    & $3t-3$ \\
$t\equiv 0$ mod $8$    & $3t/4-1$ & $\{3t,3t\}$    & $3t-3$ \\ \hline
$t\equiv 2$ mod $4$    & $t/2-1$  & $\{3t/2,3t\}$  & $3t-5$ \\
$t\equiv 12$ mod $16$  & $t/4-1$  & $\{3t/2,3t\}$  & $3t-5$ \\
$t\equiv 4$ mod $16$   & $3t/4-1$ & $\{3t/2,3t\}$  & $3t-5$ \\ \hline
$t\equiv 4$ mod $16$   & $t/4-1$  & $\{3t/4,3t\}$  & $3t-9$ \\
$t\equiv 12$ mod $16$  & $3t/4-1$ & $\{3t/4,3t\}$  & $3t-9$ \\
\end{tabular}
\caption{Orientably-regular embeddings of $K_4^{(t)}$.} \label{tab:q=4}
\end{table}

\section{Concluding remarks: chirality, duality, \\ rotational powers, and Cayley maps}\label{sec:conc}

To wrap up, we will briefly discuss the items listed in the title of this section, with motivation generated by their consideration in the original paper \cite{JaJo}.
\smallskip

By the general theory of maps \cite{JoSi}, an orientably-regular map ${\cal M}$ of type $\{m,n\}$ with orientation-preserving automorphism group $G={\rm Aut}^+{\cal M}=\langle x,y\ |\ x^m,y^n,(xy)^2,\ldots \rangle$ is {\em reflexible} (isomorphic to its `mirror image') if and only if the group admits an automorphism inverting both $x$ and $y$; in the opposite case the map is {\em chiral}. Similarly, a map ${\cal M}$ as above is {\em self-dual} if and only if its group $G={\rm Aut}^+{\cal M}$ admits an automorphism that mutually interchanges the generators $x$ and $y$. We will begin with chirality and self-duality, deferring rotational powers and Cayley maps to the second part of our concluding remarks.
\bigskip

\noindent {\bf Chirality.} The situation for orientably-regular embeddings of the graphs $K_q^{(t)}$ essentially copies the behaviour for simple complete maps from \cite{JaJo}, where it was also pointed out that chiral pairs correspond to pairs of mutually reciprocal primitive polynomials.

\begin{proposition}\label{prop:refl}
For $q\in \{2,4\}$ and any $t\ge 1$, and also for $q=3$ and any odd $t\ge 1$, all orientably-regular embeddings of $K_q^{(t)}$ are reflexible. For any even $q\ge 8$ and any $t\ge 1$, and for any odd $q\ge 5$ and odd $t\ge 1$, all orientably-regular embeddings of $K_q^{(t)}$ are chiral.
\end{proposition}

{\bf Proof.} With $G={\rm Aut}^+{\cal M}=\langle x,y\rangle$ as above, it may be checked for $q\in \{2,3,4\}$ and all admissible values of $t$ the presentations of $G$ from \eqref{eq:q=4} and \eqref{eq:Gqt} -- \eqref{eq:rel3t} are preserved by inverting $x$ and $y$, so that this inversion extends to an automorphism of $G$ and establishes reflexibility of the maps. For $q\ge 5$ and all admissible $t$, suppose that one of the orientably-regular maps from Theorem \ref{thm:Gpt} or Proposition \ref{prop:q=4} was reflexible, that is, that there was an automorphism $\theta$ of $G$ inverting both $x$ and $y$. Since the central subgroup $N=\langle y^{q-1} \rangle$ would be preserved by $\theta$, the automorphism would project onto $G/N\cong {\rm AGL}(1,q)$ and still invert the corresponding generators. This would mean that for $q\ge 5$ the regular embeddings of {\em simple} complete graphs of order $q$ were reflexible, contrary to the findings of \cite{JaJo}. \hfill $\Box$
\bigskip

\noindent {\bf Self-duality.} As it turns out, self-duality of orientably-regular embeddings of complete multigraphs behaves differently from \cite{JaJo} only for $q\le 4$.

\begin{proposition}\label{prop:dual}
An orientably-regular embedding of $K_q^{(t)}$ is self-dual if and only if the prime power $q$ and multiplicity $t$  satisfy one of the following conditions:
\medskip

\noindent {\rm (i)\ \ } $q=2$, and, in {\rm \eqref{eq:q=2}}, either $t\equiv 0$ {\rm mod} $2$ and $f=1$, or $t\equiv 0$ {\rm mod} $8$ and $f=1+t/2$;
\medskip

\noindent {\rm (ii)\ } $q=4$, and, in {\rm \eqref{eq:q=4}}, either $a=-1$ and $t\ge 1$, or $a=-1+t/2$ and $t\equiv 0$ {\rm mod} $4$, or else $a=-1\pm t/4$ and $t\equiv 8$ {\rm mod} $16$;
\medskip

\noindent{\rm (iii)} $q\in \{5,9\}$ and $t\equiv 1$ {\rm mod} $2$.
\end{proposition}

{\bf Proof}. A necessary condition for a map of type $\{m,n\}$ to be self-dual is $m=n$; this will be observed throughout in what follows, considering $q\in\{2,4\}$ first (note that $m\ne n$ for $q=3$) and then developing an argument for arbitrary $q > 4$.
\smallskip

For $q=2$, a map corresponding to a presentation \eqref{eq:q=2} can only be self-dual if $t=m=2t/\gcd(f+1,t)$, which happens only if $f=1$ (for arbitrary even $t\ge 2$) or if $f=1+t/2$ (and then $t$ must be a multiple of $4$). In the first case, \eqref{eq:q=2} reduces to $\langle x,y\ |\ x^t,\,y^t,\,(xy)^2,x^2y^2\rangle \cong C_t\times C_2$, which is clearly invariant under the interchange of $x$ and $y$. In the second case the presentation can be given the form $\langle x,y\ |\ x^t,\,y^t,\,(xy)^2,x^2y^{2-t/2}\rangle$. Raising the relation $x^2=y^{t/2-2}$ to the power of $(t/4{-}1)$ gives $x^{t/2-2}=y^2$ if $8\mid t$ and $x^{t/2-2}=y^{t/2+2}$ if $t\equiv 4$ mod $8$. One thus has invariance under the interchange of $x$ and $y$ if $t\equiv 0$ mod $8$, but not if $t\equiv 4$ mod $8$, in which case invariance would require the right-hand side of the last equation to be equal to $y^2$, contrary to the assumed order of $y$. This proves (i).
\smallskip

If $q=4$, the situation when $m=n$ refers to the presentation \eqref{eq:q=4} in which values of $a$ and $t$ come from the first four lines of Table \ref{tab:q=4} (and with $e$ as in Proposition \ref{prop:q=4}). If $a=-1$ and $e=0$, then \eqref{eq:q=4} is trivially invariant under the interchange of $x$ and $y$. If $a=-1+t/2$ for $t\equiv 0$ mod $4$, then raising the relator $x^3y^{3(1-t/2)}$ to the power of $(1{-}t/2)$ gives $y^3x^{3(1-t/2)}$ and hence implies the same invariance. Let now $a=-1+t/4$ and $t\equiv 0$ mod $8$; by a remark towards the end of the proof of Proposition \ref{prop:q=4} the relator containing $e=t/2$ can be ignored. If $t\equiv 8$ mod $16$, then raising the relator $x^3y^{3(1-t/4)}$ to the power of $(1{-}t/4)$ implies that $y^3x^{3(1-t/4)}$ and hence the above invariance. But if $t\equiv 0$ mod $16$, the same operation results in $x^{3(1-t/4)}y^{3(1-t/2)}$. In this case, if the presentation \eqref{eq:q=4} was invariant under interchange of $x$ and $y$, then there would also be the relator $y^3x^{3(1 - t/4)}$ obtained from the one we begun with. But the last two relators imply $y^{3t/2}=1$, a contradiction. The argument for $a=-1-t/4$ is entirely similar; this shows validity of (ii).
\smallskip

For $q>4$ we examine effects of pairs $(\ell,a)$ from Theorem \ref{thm:Gpt} on the relator $R=x^\ell y^{(1-q)a}$ in the presentation \eqref{eq:Gqt} -- \eqref{eq:rel3t} of $G=\langle x,y\rangle$ in the cases when valency and face length are the same, as displayed in lines 1, 2 and 4 of Table \ref{tab:q-gen}. The first and the fourth line of the table correspond to $(\ell,a)=(q-1,-1)$ for every even $q\ge 8$ and every $t\ge 1$, and also for every odd $q\equiv 1$ mod $4$ and odd $t$. The relator $R$ then becomes $x^{q-1}y^{q-1}$, which implies that the cental subgroup $N=\langle y^{q-1}\rangle< G$ is {\em equal} to $\langle x^{q-1}\rangle$. For the pair $(\ell,a)= (q-1,t/2-1)$ for even $q\ge 8$ and $t\equiv 0$ mod $4$ from the third line of Table \ref{tab:q-gen} the relator $R$ gives $x^{q-1}= y^{(q-1)(-1+t/2)}$. Since $4\mid t$, raising this relation to the power of $(t/2{-}1)$ gives $x^{(q-1)(-1+t/2)} = y^{q-1}$, implying again that $N=\langle x^{q-1}\rangle$.
\smallskip

Suppose now that an orientably-regular embedding ${\cal M}$ of $K_q^{(t)}$ for some $q>4$ and $t\ge 1$ was self-dual. This would mean that its group $G={\rm Aut}^+{\cal M}=\langle x,y\rangle$ from Theorem \ref{thm:Gpt} and $q,t$ as in lines 1,2 and 4 of Table \ref{tab:q-gen} would admit an automorphism $\theta$ interchanging $x$ and $y$. But the arguments in the previous paragraph show that in all such cases the central subgroup $N$ satisfies $N= \langle y^{q-1}\rangle = \langle x^{q-1}\rangle$, that is, $N$ is $\theta$-invariant. It follows that $\theta$ projects on to an automorphism of $G/N$ interchanging the generators $Nx$ and $Ny$, implying that the corresponding orientably-regular embeddings of {\em simple} complete graphs are self-dual. By \cite{JaJo}, however, for $q>4$ this is the case if and only if $q\in \{5,9\}$.
\smallskip

It remains to show (iii). For $q\in \{5,9\}$ let $\mu$ be a corresponding primitive polynomial; there are two choices for each $q$. Invoking Theorem \ref{thm:iso} and the notation therein but omitting subscripts, for such a pair $(q,\mu)$ and arbitrary odd $t\ge 1$ there is a {\em unique} orientably-regular map ${\cal M}(q,t;\mu)$ arising as a $t$-inflation of an orientably-regular embedding ${\cal M}(q,1;\mu)$ of a {\em simple} complete graph $K_q$, in which  $t=1$. For fixed $q\in \{5,9\}$, presentations \eqref{eq:rel1} -- \eqref{eq:rel3} of the groups $G_{q;\,\mu}= G_{q;\,\mu} (u,v) ={\rm Aut}^+{\cal M}(q,1;\mu)\cong {\rm AGL}(1,q)$ differ only in the last relator \eqref{eq:rel3} that depends on $\mu$. By Proposition \ref{prop:Schur+}, for every odd $t\ge 1$ the group $G_{q,t;\,\mu}= G_{q,t;\,\mu}(x,y) = {\rm Aut}^+{\cal M}(q,t;\,\mu)$ of the $t$-inflated map is a central extension of $G_{q;\,\mu}$ by $C_t\cong \langle y^{q-1}\rangle$. But by \cite{JaJo}, for $q\in \{5,9\}$ the original groups $G_{q;\,\mu} (u,v)$ are invariant under interchange of $u$ and $v$, that is, $G_{q;\,\mu} (u,v)\cong G_{q;\,\mu} (v,u)$. To proceed, recall our notational convention of section \ref{sec:multi-compl} by which the generators $u$ and $v$ of the original group turn into respective generators $x$ and $y$ in the extension. With this in mind, construction and uniqueness implied by Theorems \ref{thm:Gpt} and \ref{thm:iso} guarantee that the central extensions $G_{q;\,\mu} (u,v)\langle y^{q-1}\rangle = G_{q,t;\,\mu}(x,y)$ and $G_{q;\,\mu} (v,u)\langle x^{q-1}\rangle = G_{q,t;\,\mu}(y,x)$, obtained from each other by interchanging $x$ and $y$ in their presentations, determine a {\em unique} central extension of the group $G_{q;\,\mu} (u,v)\cong G_{q;\,\mu} (v,u)$. It follows that relators in the presentation \eqref{eq:Gqt} -- \eqref{eq:rel3t} determining the group $G_{q,t;\,\mu}(y,x)$ must be consequences of those obtained from its isomorphic copy $G_{q,t;\,\mu}(x,y)$ by interchanging $x$ and $y$, and vice versa. Thus, for odd $t\ge 1$ and $q\in \{5,9\}$,  $t$-inflations of orientably-regular embeddings of $K_q$ are self-dual. \hfill $\Box$
\bigskip

\noindent {\bf Rotational powers.} The group $G=G(x,y)={\rm Aut}^+{\cal M}$ of an orientably-regular map ${\cal M}$ can also be described in the form $G=G(y,z)$ in terms of generators $y$ and $z=xy$, the latter representing a `half-turn' about the centre of an edge. For any $j$ coprime to the valency $n$ of the map, one may form a new orientably-regular map ${\cal M}^j$ by replacing the $n$-fold rotation $y$ about a vertex with its $j$-th power $y^j$ while keeping $z$ intact, so that ${\rm Aut}^+{\cal M}^j= G(y^j,z)$. Because of the way of its creation the map ${\cal M}^j$ has been called the {\em $j$-th rotational power of} ${\cal M}$ in \cite{Si-surv}, but the idea of considering a `$j$-th order hole' in maps comes from \cite{CoMo}; the assignment ${\cal M}\mapsto {\cal M}^j$ is also known as the {\em Wilson operator} by its study in \cite{Wils}. The maps ${\cal M}$ and ${\cal M}^j$ have the same underlying graph and abstractly isomorphic automorphism groups, but in general they have different face length and carrier surfaces.
\smallskip

By \cite{JaJo}, for every $q\ge 2$ any two orientably-regular embeddings of $K_q$ can be obtained from each other by some rotational power. We first briefly indicate a way to extend this to $t$-inflations of such maps for every $q>2$ and every odd $t>1$. Let us begin by observing that for every $j$ with $\gcd(j,q-1) = 1$ the number $j+(q-1)\tau$ is prime for infinitely many non-negative integers $\tau$, by Dirichlet's theorem. It follows that for every $t>1$ there is a $\tau\in \{0,1,\ldots,t-1\}$ such that $\gcd(j+(q-1) \tau,(q-1)t)=1$. Let now ${\cal M}^{(t)}$ be a $t$-inflation of an orientably-regular embedding ${\cal M}$ of $K_q$ for $q>2$ and odd $t>1$. Using the description of $t$-inflations in Theorem \ref{thm:Gpt} for $q>4$ and in section \ref{sec:sing} for $q\in \{3,4\}$, one may check that the {\em uniquely determined} $(j+(q-1)\tau)$-th rotational power of ${\cal M}^{(t)}$ is isomorphic to the unique $t$-inflation of ${\cal M}^j$. Theorem \ref{thm:iso} together with the aforementioned result of \cite{JaJo} on rotational powers then imply that, for every $q>2$ and every odd $t>1$, orientably-regular embeddings of $K_q^{(t)}$ form a single orbit under the Wilson operator.
\smallskip

The `unique lifting' property may also be used for even $q\ge 8$ and even $t\ge 2$, to be outlined next. For such $q$ and $t$ the key is to use factorisation of the automorphism group by the central subgroup $\langle y^{2(q-1)}\rangle$, that is, by `half' of $N=\langle y^{q-1}\rangle$ used throughout. By Theorem \ref{thm:iso}, for $t=2$ and $q=2^k$ with $k\ge 3$ this results in $2\varphi(q-1)/k$ pairwise non-isomorphic maps ${\cal M}_c(q,2;\mu)$, distinguished by a choice of one of the $\varphi(q-1)/k$ primitive polynomials $\mu$ and a choice of the value of the parameter $c\in \{0,1\}$ from part (II) of Theorem \ref{thm:Gpt}; each such maps lifts to a {\em unique} map ${\cal M}_c(q,t;\mu)$. In terms of $y$ and $z$, presentations {\rm \eqref{eq:Gqt} -- \eqref{eq:rel3t}} of the groups ${\rm Aut}^+{\cal M}_c(q,2;\mu)$ for $t=2$, $a=c-1$, $1\le i\le k-1$, and with $W(y,z;\mu)$ obtained from $W(x,y;\mu)$ by letting $x=zy^{-1}$, reduce to
\begin{equation}\label{eq:t=2}
\langle y,\,z\ |\ y^{2(q-1)},z^2,(yz)^{q-1}y^{(c-1)(q-1)},[y^{q-1},z],\ldots,[z,y^izy^{-i}],\ldots,W(y,z;\mu)y^{(1{-}q)c} \rangle
\end{equation}

Let $j$ be an integer coprime to $2(q-1)$, or, equivalently, and odd integer coprime to $q-1$. By the background on rotational powers, for $c\in \{0,1\}$ the automorphism group of a $j$-th rotational power of one of the maps ${\cal M}_c(q,2;\mu)$ is obtained by substituting $y^j$ for $y$ in the presentation \eqref{eq:t=2}. We will focus on the effect of such a substitution in the third relator of \eqref{eq:t=2} by evaluating $(y^jz)^{q-1}$ in terms of $(yz)^{q-1}$. To do so we will use the commutation relators $[z,y^izy^{-i}]$ which, by a remark in the last paragraph of the proof of Theorem \ref{thm:pres} and by $y^{q-1}$ commuting with $z$, extend to relators $[y^{i'}zy^{-i'}, y^izy^{-i}]$ for any integers $i,i'$. Letting $z_i=y^izy^{-i}$ for arbitrary $i$, one has $(y^jz)^{q-1}=z_jz_{2j}\dots z_{(q-1)j}y^{(q-1)j}$. But the sets $\{j,\,2j,\,\ldots,\,(q-1)j\}$ mod $q-1$ and $\{0,\,1,\,2,\,\ldots,q-2\}$ mod $q-1$ are identical, and note that $z_{i+i'(q-1)}=z_i$ for any $i,i'$. Now, since any two of the $q$ terms on the right-hand-side of the last product commute and $y^{(q-1)j}=y^{q-1}$ by oddness of $j$ and ${\rm ord}(y)=2(q-1)$, we may rewrite the product in the form \[(y^jz)^{q-1} = z_jz_{2j}\dots z_{(q-1)j}\cdot y^{(q-1)j} = zz_{1}z_{2}\ldots z_{q-2}\cdot y^{q-1} = (zy)^{q-1}\]
and hence $(y^jz)^{q-1} = (zy)^{q-1}$. In view of the form of the third relator in \eqref{eq:t=2} this implies that any rotational power of a map of the form ${\cal M}_c(q,2;\mu)$ results in a map {\em with the same value of $c$}. Combined with the findings of \cite{JaJo} for even $q>4$ it follows that for $t=2$ and any fixed even $q>4$ there are two orbits of orientably-regular embeddings of $K_q^{(2)}$ under rotational powers, one for each $c\in \{0,1\}$. Applying the unique lifting described before for odd $t$ and arbitrary $q>2$ it follows that, for every fixed even $q>4$ and even $t\ge 2$, orientably-regular embeddings of $K_q^{(t)}$ form two orbits under the Wilson operator.
\smallskip

One may check that the conclusion from the previous paragraph is valid also for $q=4$ and $t\equiv 2$ mod $4$. The case $q=4$ and $4\mid t$ forms an exception again and can be resolved using factorisation by the central subgroup $\langle y^{4(q-1)} \rangle$. This results in {\em three} orbits of regular embeddings of $K_4^{(4)}$ for $q=t=4$: a singleton orbit for each $a\in\{-1,1\}$ from part (A) of Proposition \ref{prop:q=4}, and one orbit with two maps for $a\in \{0,-2\}$ from part (B) of Proposition \ref{prop:q=4}, arising from each other by taking the $5$th and also the $7$th rotational power. The rest follows by unique lifting as in the two cases considered above, giving one orbit for each $a\in \{-1,-1+t/2\}$ and a third orbit fusing the maps for $a=-1\pm t/4$ by $j$-th rotational powers for any $j\equiv 5$ or $7$ mod $12$ such that $\gcd(j,3t)=1$. Summing up, one obtains:

\begin{proposition}\label{prop:Wil}
For every $q>2$ and $t\ge 1$, under the Wilson operator the orientably-regular embeddings of $K_{q}^{(t)}$ form
\smallskip

\noindent{\rm (a)} a single orbit if $t$ is odd;
\smallskip

\noindent{\rm (b)} two orbits if both $q$ and $t$ are even and $q>4$, or if $q=4$ and $t\equiv 2$ mod $4$;
\smallskip

\noindent{\rm (c)} three orbits if $q=4$ and $t\equiv 0$ mod $4$.  \hfill $\Box$
\end{proposition}

For completeness we note that, for $q=2$ and every $t\ge 2$, every orientably-regular embedding of $K_2^{(t)}$ forms a singleton orbit under the Wilson operator.

\bigskip

\noindent{\bf Cayley maps.} We conclude by outlining another way of approaching inflations of embeddings of complete graphs $K_q$ for prime powers $q=p^k$, based on structural results of the seminal paper \cite{JaJo}. These imply that every orientably-regular embedding ${\cal M}$ of $K_q$ has the property that the commutator subgroup $G'\cong C_p^k\cong F_q^+$ of the group $G={\rm Aut}^+{\cal M}\cong {\rm AGL}(1,q)\cong F_q^+\rtimes F_q^{\times}$ acts regularly on the vertex set of the embedded complete graph. This not only means that one can identify the vertex set of ${\cal M}$ with $F_q^+$ via this regular action, which we will assume in what follows, but also that the (say, anti-clockwise) cyclic order of neighbours of every vertex $\nu\in F_q^+$ on the carrier surface of ${\cal M}$ has the form $(\nu+1,\nu+\xi,\nu+\xi^2,\ldots,\nu+\xi^{q-2})$ for some primitive element $\xi\in F_q$ (cf. \cite[section 5] {JaJo}). An orientably-regular map admitting a regular action of a {\em normal} subgroup of its orientation-preserving automorphism group on vertices has been known as a {\em balanced Cayley map} \cite{SkSi}. One of the outcomes of the analysis in \cite{JaJo} is that all the orientably-regular embeddings of complete graphs are balanced Cayley maps.
\smallskip

In orientably-regular maps arising as inflations of complete maps, a regular action of the corresponding commutator subgroup on vertices of an embedding ${\cal M}(q,t;\mu)$ of a complete multigraph $K_q^{(t)}$ for $q\ne 4$ is a consequence of Proposition \ref{prop:Schur+}. The vertex set of ${\cal M}(q,t;\mu)$ may in this case be identified with $F_q^+$ again, now regarded as a $k$-dimensional vector space over $F_p$ formed by column vectors as in section \ref{sec:simple}. The fact that for $t>1$ the inflated map has multiple edges turns out to be immaterial, and for every vertex ${\mathbf b}\in F_q^+$ the cyclic order of its neighbours, including their repeated occurrence, can be described via the action $C_{(q-1)t}\cong \langle y\rangle\to {\rm Aut}(F_q^+)$, $j\mapsto A^j$ from Proposition \ref{prop:Schur+} in the form $({\mathbf b}{+}{\mathbf e},\, {\mathbf b}{+}A{\mathbf e},\,\ldots, {\mathbf b}{+}A^j{\mathbf e},\,\ldots, {\mathbf b}{+}A^{t(q-1)-1}{\mathbf e})$, where ${\mathbf e}\in F_q^+{\setminus} \{ {\mathbf 0} \}$ is an arbitrary but fixed vector. Thus, for $q\ne 4$, orientably-regular embeddings of complete multigraphs on $q$ vertices are also (balanced) Cayley maps, albeit with multiple edges. This remains valid also for $q=4$ if $G'\cong C_2\times C_2$ but {\em not} if $G'$ is isomorphic to the quaternion group. Indeed, calculations in the proof of part (B) of Proposition \ref{prop:q=4} imply that in this case the map automorphism group contains no regular subgroup of order $4$ on the vertex set, underscoring a remarkable exceptionality of the situation for $q=4$ in yet another way.
\smallskip

We note that a general theory of Cayley maps, including those in which underlying graphs have multiple edges, was  developed in \cite{R+}.
\smallskip

Our last remark addresses potential demand for concrete examples. Although we did not include any, we worked out a number of particular cases in proofs of Propositions \ref{prop:q=4} -- \ref{prop:Wil}. In the light of these, additional examples of substituting concrete values of parameters in the presentation \eqref{eq:Gqt} -- \eqref{eq:rel3t} are not likely to illuminate the subject any further.
\bigskip

\bigskip

{\bf Acknowledgment.} The authors acknowledge support of this research from APVV Research Grants 19-0308 and 22-0005, and VEGA Research Grants 1/0567/22 and 1/0069/23.

\end{document}